\documentclass[12pt,twoside,leqno]{article}
\usepackage{amsmath}
\usepackage{amssymb}
\usepackage{amsxtra}
\usepackage{amscd}
\usepackage{amsthm}
\usepackage[mathscr]{eucal}
\usepackage{xr}
\theoremstyle{plain}
\newtheorem{thm}[subsection]{Theorem}
\newtheorem{prop}[subsection]{Proposition}

\newtheorem{lem}[subsection]{Lemma}
\newtheorem{conj}[subsection]{Conjecture}

\theoremstyle{definition}

\newtheorem{rem}[subsection]{Remark}
\newtheorem{para}[subsection]{}

\newenvironment{pf}{\proof[\proofname]}{\endproof}

\begin{document}
\title{Mixed objects are embedded into log pure objects}

\author{Kazuya Kato, Chikara Nakayama, and Sampei Usui}

\maketitle

\newcommand\Cal{\mathcal}
\newcommand\define{\newcommand}

\define\gp{\mathrm{gp}}%
\define\fs{\mathrm{fs}}%
\define\an{\mathrm{an}}%
\define\mult{\mathrm{mult}}%
\define\Ker{\mathrm{Ker}\,}%
\define\Coker{\mathrm{Coker}\,}%
\define\Hom{\mathrm{Hom}\,}%
\define\Ext{\mathrm{Ext}\,}%
\define\rank{\mathrm{rank}\,}%
\define\gr{\mathrm{gr}}%
\define\cHom{\Cal Hom\,}%
\define\cExt{\Cal Ext\,}%
\define\cA{\Cal A}
\define\cB{\Cal B}
\define\cC{\Cal C}
\define\cD{\Cal D}
\define\cE{\Cal E}
\define\cF{\Cal F}
\define\cG{\Cal G}
\define\cH{\Cal H}
\define\cM{\Cal M}
\define\cN{\Cal N}
\define\cO{\Cal O}
\define\cS{\Cal S}
\define\cT{\Cal T}
\define\fF{\frak F}
\define\Cc{C^{\spcheck}}
\define\Dc{\check{D}}

\newcommand{\N}{{\mathbb{N}}}
\newcommand{\Q}{{\mathbb{Q}}}
\newcommand{\Z}{{\mathbb{Z}}}
\newcommand{\R}{{\mathbb{R}}}
\newcommand{\C}{{\mathbb{C}}}
\newcommand{\bP}{{\mathbb{P}}}
\newcommand{\bN}{{\mathbb{N}}}
\newcommand{\bQ}{{\mathbb{Q}}}
\newcommand{\bF}{{\mathbb{F}}}
\newcommand{\bZ}{{\mathbb{Z}}}
\newcommand{\bR}{{\mathbb{R}}}
\newcommand{\bC}{{\mathbb{C}}}
\newcommand{\bbQ}{{\bar \mathbb{Q}}}
\newcommand{\ol}[1]{\overline{#1}}
\newcommand{\too}{\longrightarrow}
\newcommand{\respect}{\rightsquigarrow}
\newcommand{\compatible}{\leftrightsquigarrow}
\newcommand{\upc}[1]{\overset {\lower 0.3ex \hbox{${\;}_{\circ}$}}{#1}}
\newcommand{\Gmlog}{\mathbb{G}_{m,\log}}
\newcommand{\Gm}{\mathbb{G}_m}
\newcommand{\ep}{\varepsilon}
\newcommand{\Spec}{\operatorname{Spec}}
\newcommand{\val}{{\mathrm{val}}} 
\newcommand{\n}{\operatorname{naive}}
\newcommand{\bs}{\operatorname{\backslash}}
\newcommand{\Gal}{\operatorname{{Gal}}}
\newcommand{\gal}{{\rm {Gal}}({\bar \Q}/{\Q})}
\newcommand{\galp}{{\rm {Gal}}({\bar \Q}_p/{\Q}_p)}
\newcommand{\gall}{{\rm{Gal}}({\bar \Q}_\ell/\Q_\ell)}
\newcommand{\wep}{W({\bar \Q}_p/\Q_p)}
\newcommand{\wel}{W({\bar \Q}_\ell/\Q_\ell)}
\newcommand{\Ad}{{\rm{Ad}}}
\newcommand{\BS}{{\rm {BS}}}
\newcommand{\even}{\operatorname{even}}
\newcommand{\End}{{\rm {End}}}
\newcommand{\odd}{\operatorname{odd}}
\newcommand{\GL}{\operatorname{GL}}
\newcommand{\Lam}{{\Lambda}}
\newcommand{\La}{{\Lambda}}
\newcommand{\la}{{\lambda}}
\newcommand{\tLam}{{\tilde \Lambda}}
\newcommand{\tLa}{{\tilde \Lambda}}
\newcommand{\uL}{{{\hat {L}}^{\rm {ur}}}}
\newcommand{\uQp}{{{\hat \Q}_p}^{\text{ur}}}
\newcommand{\sel}{\operatorname{Sel}}
\newcommand{\dt}{{\rm{Det}}}
\newcommand{\Sig}{\Sigma}
\newcommand{\dla}{{\rm{Det}}_{\Lambda}}
\newcommand{\can}{{\rm {can}}}
\newcommand{\dSla}{{\rm{Det}}_{S^{-1}\La}}
\newcommand{\tSla}{S^{-1}\La\otimes_{\La}}
\newcommand{\Sla}{S^{-1}\La}
\newcommand{\SL}{{\rm{SL}}}
\newcommand{\spl}{{\rm{spl}}}
\newcommand{\Isom}{{\rm {Isom}}}
\newcommand{\Mor}{{\rm {Mor}}}
\newcommand{\bg}{\bar{g}}
\newcommand{\id}{{\rm {id}}}
\newcommand{\cone}{{\rm {cone}}}
\newcommand{\al}{a}
\newcommand{\ChL}{{\cal{C}}(\La)}
\newcommand{\Image}{{\rm {Image}}}
\newcommand{\Dp}{{D_{\text{parf}}}}
\newcommand{\Aut}{{\rm {Aut}}}
\newcommand{\Qp}{{\mathbb{Q}}_p}
\newcommand{\barQp}{{\mathbb{Q}}_p}
\newcommand{\Qpur}{{\mathbb{Q}}_p^{\rm {ur}}}
\newcommand{\Qppin}{{\mathbb{Q}}_p(\mu_{p^{\infty}})}
\newcommand{\Zp}{{\mathbb{Z}}_p}
\newcommand{\Zl}{{\mathbb{Z}}_l}
\newcommand{\Ql}{{\mathbb{Q}}_l}
\newcommand{\Qlur}{{\mathbb{Q}}_l^{\rm {ur}}}
\newcommand{\F}{{\mathbb{F}}}
\newcommand{\eps}{{\epsilon}}
\newcommand{\epsLa}{{\epsilon}_{\La}}
\newcommand{\epsLaVxi}{{\epsilon}_{\La}(V, \xi)}
\newcommand{\epsOLaVxi}{{\epsilon}_{0,\La}(V, \xi)}
\newcommand{\Ga}{{\Gamma}}
\newcommand{\Qplin}{{\mathbb{Q}}_p(\mu_{l^{\infty}})}
\newcommand{\otimesQplin}{\otimes_{\Qp}{\mathbb{Q}}_p(\mu_{l^{\infty}})}
\newcommand{\galFl}{{\rm{Gal}}({\bar {\Bbb F}}_\ell/{\Bbb F}_\ell)}
\newcommand{\gallur}{{\rm{Gal}}({\bar \Q}_\ell/\Q_\ell^{\rm {ur}})}
\newcommand{\galFF}{{\rm {Gal}}(F_{\infty}/F)}
\newcommand{\galFv}{{\rm {Gal}}(\bar{F}_v/F_v)}
\newcommand{\galF}{{\rm {Gal}}(\bar{F}/F)}
\newcommand{\epsVxi}{{\epsilon}(V, \xi)}
\newcommand{\epsOVxi}{{\epsilon}_0(V, \xi)}
\newcommand{\ccy}{\chi_{\rm {cyclo}}}
\newcommand{\plim}{\lim_
{\scriptstyle 
\longleftarrow \atop \scriptstyle n}}
\newcommand{\sig}{{\sigma}}
\newcommand{\ga}{{\gamma}}
\newcommand{\del}{{\delta}}
\newcommand{\Vss}{V^{\rm {ss}}}
\newcommand{\Bst}{B_{\rm {st}}}
\newcommand{\Dpst}{D_{\rm {pst}}}
\newcommand{\Dcrys}{D_{\rm {crys}}}
\newcommand{\DdR}{D_{\rm {dR}}}
\newcommand{\Fin}{F_{\infty}}
\newcommand{\Kla}{K_{\lambda}}
\newcommand{\Ola}{O_{\lambda}}
\newcommand{\Mla}{M_{\lambda}}
\newcommand{\Det}{{\rm{Det}}}
\newcommand{\LaSa}{{\La_{S^*}}}
\newcommand{\cX}{{\cal {X}}}
\newcommand{\MHG}{{\frak {M}}_H(G)}
\newcommand{\tauMla}{\tau(M_{\lambda})}
\newcommand{\Fvur}{{F_v^{\rm {ur}}}}
\newcommand{\Lie}{{\rm {Lie}}}
\newcommand{\cL}{{\cal {L}}}
\newcommand{\fq}{{\frak {q}}}
\newcommand{\cont}{{\rm {cont}}}
\newcommand{\SC}{{SC}}
\newcommand{\Om}{{\Omega}}
\newcommand{\dR}{{\rm {dR}}}
\newcommand{\hatSig}{{\hat{\Sigma}}}
\newcommand{\Ups}{{\Upsilon}}
\newcommand{\bUps}{{{\bar{\Upsilon}}}}
\newcommand{\rdet}{{{\rm {det}}}}
\newcommand{\ord}{{{\rm {ord}}}}
\newcommand{\BdR}{{B_{\rm {dR}}}}
\newcommand{\BdRO}{{B^0_{\rm {dR}}}}
\newcommand{\Bcrys}{{B_{\rm {crys}}}}
\newcommand{\Qw}{{\mathbb{Q}}_w}
\newcommand{\barkappa}{{\bar{\kappa}}}
\newcommand{\cP}{{\Cal {P}}}
\newcommand{\oppLa}{{\Lambda^{\circ}}}
\newcommand{\LMH}{{\rm {LMH}}}
\newcommand{\spe}{{\rm {sp}}}
\newcommand{\BM}{{\rm {BM}}}
\newcommand{\pts}{{\rm {pts}}}
\newcommand{\pe}{\frak p}
\newcommand{\adm}{{\rm {adm}}}
\newcommand{\triv}{{\rm {triv}}}
\newcommand{\et}{{\rm {\acute et}}}

\newcommand{\LM}{{\rm {LM}}}
\newcommand{\MM}{{\rm {MM}}}
\newcommand{\MH}{{\rm {MH}}}
\newcommand{\LH}{{\rm {LH}}}

\begin{abstract}
\noindent We prove that a variation of mixed Hodge structure is embedded in a logarithmic variation of pure Hodge structure, and a generalized version of this result. 
  These results suggest some simple construction 
of the category of mixed motives by using log pure motives. 
\end{abstract}
\renewcommand{\thefootnote}{\fnsymbol{footnote}}
\footnote[0]{Primary 14A21; Secondary 14D07, 32G20}

\section*{Introduction}
\begin{para}  In this paper, we develop our idea in \cite{KNU08} that a mixed 
object is embedded in  a log pure object. We improve the result in \cite{KNU08} on this idea (Theorem \ref{MP}) and propose a simple construction  of the category of mixed motives over a field based on this idea without assuming any conjecture (Appendix in this paper). 
  
\end{para}

\begin{para}
The following is a standard example concerning this idea.

 Let $\Delta$ be the unit disc $\{q\in \C\;|\; |q|<1\}$, and let $\frak X$ be a smooth complex manifold with a projective flat morphism $\frak X\to \Delta$  which is smooth outside $0\in \Delta$ and is of semistable reduction at $0\in \Delta$. For $t\in \Delta$, let  $\frak X_t\subset \frak X$ be the fiber over $t\in \Delta$. Then we have the mixed Hodge structure $H^1(\frak X_0, \Z)$. This mixed Hodge structure is embedded in the limit mixed Hodge structure $``\lim_{t\neq 0, t\to 0}\;  H^1(\frak X_t,\Z)"\supset  H^1(\frak X_0, \Z)$, and this limit mixed Hodge structure is associated with  
the  log pure Hodge structure  $H=H^1((\frak X_0\;\text{with log})/(0 \;\text{with log}),\Z)$ of weight $1$ on the standard log point $0\in \Delta$. Thus the mixed object $H^1(\frak X_0,\Z)$  is embedded in the  log pure object $H$. 

\end{para}

\begin{para} In \cite{KNU08}, we proved that a mixed Hodge structure is embedded in a log pure Hodge structure, which is the case $n=0$ of the following more general result proved in \cite{KNU08}: A nilpotent orbit of mixed Hodge structures with $n$ monodromy operators is embedded in a nilpotent orbit of pure Hodge structures with one more monodromy operators. This general result was successfully applied in \cite{KNU08} to deduce the $\SL(2)$-orbit theorem for the degeneration of mixed Hodge structure from the $\SL(2)$-orbit theorem of Cattani--Kaplan--Schmid (\cite{CKS1}) for degeneration of pure Hodge structure.

In this paper, we prove the following  further 
generalization (Theorem \ref{MP} in Section 1) of the  result in \cite{KNU08}:  A log mixed Hodge structure on an fs log analytic space $X$ with polarizable graded quotients for the weight filtration is, locally on $X$,  embedded into a log pure Hodge structure on $X\times S$, where $S$ is the standard log point. The text of this paper (Sections \ref{s:open}--\ref{sec:pf}) is  devoted to the proof of this theorem.

\end{para}

\begin{para}

In the theory of mixed motives over a field $k$, a big question is how to define the set of morphisms of mixed motives
$$(1) \quad h(Y)(r) \to h(Z)(s)$$
for schemes $Y, Z$ of finite type over $k$ and for $r,s\in \Z$.  Here $h(Y)$ is the mixed motive associated with $Y$ whose $\ell$-adic realization for a prime number $\ell\neq \text{char}(k)$ is $\bigoplus_m\; H^m_{\et}(Y\otimes_k \overline k, \Q_{\ell})$ and $(r)$ means the Tate twist. By the construction of the category (\MM$**$)  of mixed motives over $k$ in \ref{**} in Appendix, we answer this question as follows. We define the category $(\LM\flat\text{$**$})$ of limit mixed motives associated with log pure motives by using certain $K$-groups as the sets of morphisms, and 
define the mixed motive $h(Y)(r)$ as a functor from $(\LM\flat\text{$**$})$ to the category of $\Q$-vector spaces, by using certain $K$-groups. Thus a morphism (1) is a morphism of functors. 

We hope that this method is justified by its Hodge version (\ref{Hodge}): Our result on Hodge theory 
tells in particular that we can regard a mixed Hodge structure  as a functor on the category of limit mixed Hodge structures associated with log pure  Hodge structures.

The notion mixed motive is more difficult than the notion log pure motive (the latter is just the logarithmic version of the pure motive of Grothendieck) and our hope is that the difficult objects  mixed motives are well-understood by using log pure motives which are simpler.

  This Appendix (Section A), which discusses the motive theory, is independent of the text and one can read it first. 
  
\end{para}

  K.\ Kato was 
partially supported by NFS grants DMS 1303421, DMS 1601861, and DMS 2001182.
C.\ Nakayama was 
partially supported by JSPS Grants-in-Aid for Scientific Research (B) 23340008, (C) 16K05093, and (C) 21K03199.
S.\ Usui was 
partially supported by JSPS Grants-in-Aid for Scientific Research (B) 23340008, 
(C) 17K05200, and (C) 22K03247. 

\bigskip

\section{The results}\label{sec:main}

\begin{para}\label{MP0}  
As in \cite{KU}, let $\cA(\log)$ be the category of  fs log analytic spaces (i.e., complex analytic spaces with fs log structures) and let $\cB(\log)\supset \cA(\log)$ be the category of locally ringed spaces over $\C$ with log structures which are locally subspaces of objects of $\cA(\log)$ with the strong topologies (\cite{KU} 3.2). 

Fix a subring $R$ of $\R$.

Let $X$ be an object of $\cB(\log)$ and let $H$ be as in one of the following (1) and (2).

\medskip

(1) $H$ is an $R$-log mixed Hodge structure ($R$-LMH) on $X$.

(2) $H$ is an $R$-log variation of mixed Hodge structure ($R$-LVMH) on $X$. 

\medskip
For the definitions of $R$-LMH and $R$-LVMH, for the definitions of $R$-polarized log Hodge structure ($R$-PLH) and $R$-log variation of polarized Hodge structure ($R$-LVPH), and for the pre-versions (pre-$R$-LMH, etc.), cf.\ \cite{KU} 2.6 and \cite{KNU3} 1.3, where the cases $R=\bZ$ are treated.
 The difference of $R$-LMH (resp.\ $R$-PLH) and $R$-LVMH (resp.\ $R$-LVPH) lies in that the latter must satisfy the big Griffiths transversality though the small Griffiths transversality is satisfied by the former (\cite{KU} 2.4.9).

  In both situations (1) and (2), we assume that $H$ satisfies the following conditions (i) and (ii).

\medskip

(i)  The local system $H_R$ is locally free as a sheaf of $R$-modules on $X^{\log}$. Furthermore, $W_wH_R:= H_R\cap W_w(H_R \otimes_{\Z} \Q)\subset H_R \otimes_{\Z} \Q$  and $\gr^W_wH_R:= W_wH_R/W_{w-1}H_R$ for all $w$  are locally free as sheaves of $R$-modules on $X^{\log}$. 

(ii) For each $w$, there is an $R$-perfect $(-1)^w$-symmetric bilinear form $\gr^W_wH_R \times \gr^W_wH_R \to R\cdot(2\pi i)^{-w}$ which gives a polarization of $\gr^W_wH$.

\medskip

 If $R$ is a field (as in the important cases $R=\Q$, $R=\R$), the condition  (i) is empty  and  the condition (ii) simply says that $H$ has polarizable $\gr^W$. 

\end{para}

The aim of this paper is to prove 
\
\begin{thm}\label{MP}

  Assume that we are in the situation $(1)$ (resp.\ $(2)$) in $\ref{MP0}$. 
  Let $S$ be the standard log point. 
  Then locally on $X$, there are an $R$-PLH (resp.\ $R$-LVPH)
$H'$ on $X\times S$ and an injective homomorphism $H_R \to H'_R$ of the local systems of $R$-modules on $(X\times S)^{\log}$ satisfying the following conditions {\rm (i)}, {\rm (ii)}, {\rm (a)}, and {\rm (b)} below. If $W_wH=H$, there is such an $H'$ of weight $w$.

\medskip

{\rm (i)}  $ H'_R$ and $H'_R/H_R$ are locally free as sheaves of  $R$-modules on $(X\times S)^{\log}$. 

{\rm (ii)}  The polarization of $H'$ is given by an $R$-perfect $(-1)^{w'}$-symmetric bilinear form $H'_R\times  H'_R\to R\cdot(2\pi i)^{-w'}$, where $w'$ is the weight of $H'$. 

(These conditions {\rm (i)} and {\rm (ii)} are automatically satisfied if $R$ is a field.)

\medskip

{\rm (a)} The Hodge filtration of $H$ is the restriction of that of $H'$. More precisely, the Hodge filtration of $H$ on $H_{\cO}= (\tau_X)_*(\cO^{\log}_X\otimes_R H_R)= (\tau_{X\times S})_*(\cO^{\log}_{X\times S} \otimes_R H_R)$ coincides with the restriction of the Hodge filtration of $H'$ on $H'_{\cO}=(\tau_{X\times S})_*(\cO^{\log}_{X\times S} \otimes_R H'_R)$.

{\rm (b)} The weight filtration of $H$ is the restriction of the relative monodromy filtration of $H'$. More precisely, for every $t\in (X\times S)^{\log}$,  the  weight filtration of $H$ on the stalk $H_{R, t}\otimes_\Z \Q$  is the restriction of the relative monodromy filtration on $H'_{R,t}\otimes_{\Z} \Q$ of the logarithm $H'_{R, t}\otimes_{\Z} \Q \to  H'_{R,t}\otimes_\Z \Q$ of the action of the standard generator of 
$\pi_1(S^{\log})$. 

\end{thm}

\begin{rem}\label{MPrem}

(1) By duality, we have a result in which we replace the injection $H_R\to  H'_R$ in Theorem \ref{MP} by a surjection $H'_R\to H_R$ and change the conditions (i), (a), and (b) accordingly.

(2)  \cite{KNU08} Proposition 4.1 is a slightly weaker version of the case where $X$ is an fs log point 
of this theorem. The structure of the proof of the above theorem given below is similar to that of the proof of \cite{KNU08} Proposition 4.1 given in Section 6 and Section 7 of \cite{KNU08}.

(3)  On the other hand, the case of Theorem \ref{MP} where $X=(\Spec\bC, \bC\oplus \bN^n)$ implies that we can take all the $a_{jk}$ to be $0$ unless $j=k$ in \cite{KNU08} Proposition 4.1 (cf.\ the remark after ibid.\ Proposition 4.1). 
  As explained in ibid.\ 5.9, this gives a characterization of 
$\bR$-IMHM (\cite{KNU08} 5.2, \cite{Kas86}) without using relative monodromy filtrations. 
  We state this below as Proposition \ref{p:IMHM}. 

(4) When $X$ has the trivial log structure, this theorem implies the following. 
A variation of mixed Hodge structure with polarizable graded quotients 
on a 
complex analytic manifold 
$X$ is, locally on $X$, embedded in  a 
log variation of polarized Hodge structure. 
\end{rem}

\begin{prop}
\label{p:IMHM}
  Let $(V, W, N_1,\ldots, N_n, F)$ be a pre-$\bR$-IMHM ({\rm \cite{KNU08}} $5.2$). 
  It is an $\bR$-IMHM if and only if there is a pure nilpotent orbit 
$(V', w, N'_0,\ldots, N'_n, F')$ and a surjective homomorphism 
$(V', W(N'_0)[-w], N'_1,\ldots, N'_n, F') \to (V, W, N_1,\ldots, N_n, F)$ 
of pre-$\bR$-IMHMs.
\end{prop}

\begin{para}\label{insp2} Inspired by Remark \ref{MPrem} (4), we  expect that a motive theoretic version of the above theorem exists, that is, 
that a mixed motive can be embedded into a log pure motive. Based on this idea, we  construct
 the category of mixed motives over a field in Appendix (Section \ref{s:motive}) by using log pure motives. 

\end{para}

\section{Preparation on log Hodge theory}
\label{s:open}
We prove two propositions on log Hodge theory together. Proposition \ref{open} will be used in the last part of the proof of Theorem \ref{Pmixed0}. Proposition \ref{fromLM} will be used in the proof of Lemma \ref{dual}.

\begin{prop}\label{open} Let $X$ be an object of $\cB(\log)$, let $R$ be a subfield of $\R$,  and let $H$ be a  pre-$R$-LMH  on $X$ satisfying the small Griffiths transversality. Assume that for each $w\in \Z$, we are given a $(-1)^w$-symmetric pairing $\langle\cdot,\cdot\rangle_w: \gr^W_wH \otimes \gr^W_wH\to R(-w)$ which induces an isomorphism $\gr^W_wH\overset{\cong}\to(\gr^W_w  H)^*(-w)$ of pre-$R$-LMH, where $(\cdot)^*$ denotes the dual. Let $U$ be the set of all $x\in X$ such that the pullback of $H$ to the fs log point $x$ is an  $R$-LMH and such that $\langle\cdot,\cdot\rangle_w$ is a polarization for every $w\in \Z$. Then $U$ is an open  set of  $X$.

\end{prop}

\begin{prop}\label{fromLM}  Let $X$ be an object of $\cB(\log)$, let $R$ be a subring of $\R$, and let $H$ be an $R$-LMH  on $X$ 
with polarized $\gr^W$ satisfying the condition $\ref{MP0}$ {\rm (i)}.
Then locally on $X$, there are a log manifold $Z$ (for a log manifold, see {\rm \cite{KU}} Definition {\rm $3.5.7$}) and a morphism $X\to Z$ of $\cB(\log)$ 
 such that $H$ is the pullback of an $R$-LMH on $Z$
 with  polarized  $\gr^W$  satisfying the condition $\ref{MP0}$  {\rm (i)}.

\end{prop}

\begin{rem}

In the case of $R=\Z$ or $\Q$, under the assumption that the local monodromy of $H$ at each point of $X$ is contained in a sharp cone, Proposition \ref{fromLM}  is a consequence of the existence of the moduli space of LMH with polarized $\gr^W$ and of the fact that this moduli space is a log manifold as treated in \cite{KU} and \cite{KNU3}. Here we treat an $R$-LMH without such an assumption on local monodromy. 
The proof of Proposition \ref{fromLM} uses arguments in \cite{KU} 2.3.7 and \cite{KU} Section 8, and the space $E$ (resp.\ $\check{E}$) which appears in \ref{E1} below is a variant 
of the space  $E_{\sig}$ (resp.\ $\check{E}_{\sig}$) in \cite{KU} and \cite{KNU3}.

\end{rem}

\begin{para}\label{pfLM} Let $X$ be an object of $\cB(\log)$. Assume that we are given a pre-$R$-LMH $H$ on $X$. For the proof of Proposition \ref{open} (resp.\ \ref{fromLM}), we assume that we are given  $\langle\cdot, \cdot\rangle_w$ on $\gr^W_wH$ for each $w$ which is as in the hypothesis of Proposition \ref{open} (resp.\ which is a polarization).

\end{para}

\begin{para}\label{xi1} Let $s\in X$ and let $t$ be a point of $X^{\log}$ lying over $s$. 
  We work around $s$. 
  Let $(q_j)_{1\leq j\leq n}$ be a finite family of local sections of $M_X$  around $s$ which forms a $\Z$-base of $(M_X^{\gp}/\cO_X^\times)_s$. For $1\leq j\leq n$, let $\gamma_j$ be the element of $\pi_1(s^{\log})=\Hom((M_X^{\gp}/\cO_X^\times)_s,\Z)$ which sends $q_j$ to $1$ and $q_k$ to $0$ for all $k\neq j$
(see \cite{KU} 2.2.9 for this identification).
Then the  action of $\gamma_j$ on $H_{R,t}$ is unipotent.  Let $N_j=(\log(\gamma_j)): H_{R, t}\otimes_{\Z}\Q\to H_{R,t}\otimes_{\Z} \Q$. Let $z_j$ be a branch of $ (2\pi i)^{-1}\log(q_j)$ around  $t$. Then on an open neighborhood of $t$ in $X^{\log}$, concerning $H_{\cO}=\tau_*(\cO^{\log}_X \otimes_R H_R)$, we have 
$$H_{\cO}= \exp\bigl(\sum_{j=1}^n z_jN_j\bigr)(\cO_X\otimes_R H_R).$$

\end{para}

\begin{para}\label{cS}
  By replacing $X$ by an open neighborhood of $s$ in $X$ if necessary, we may assume that 
there is  a chart $\cS\to M_X$ with an fs monoid $\cS$ such that $\cS\overset{\cong}\to (M_X/\cO^\times_X)_s$. 
Let $\cT=\Spec(\C[\cS])^{\an}$ and let $X\to \cT$ be the morphism induced by the composition $\cS\to M_X\to \cO_X$. This induces an isomorphism from $s$ to the ``origin'' of $\cT$. Since the induced map $s^{\log}\to \cT^{\log}$ is a homotopy equivalence, the restriction of $H_R$ to $s^{\log}$ extends uniquely to a local system $H_{R, \cT}$ on $\cT^{\log}$ which also has a weight filtration $W$ and the family $(\langle\cdot, \cdot\rangle_w)_w$ of pairings. By the properness of $X^{\log}\to X$, for some open neighborhood $V$ of $s$ in $X$, we have an isomorphism between the pullbacks of $(H_R,W, (\langle\cdot,\cdot\rangle_w)_w)$ and  $(H_{R,\cT}, W, (\langle\cdot,\cdot\rangle_w)_w)$ to $V^{\log}$.

\end{para}

\begin{para} 
Let $H_0:=H_{R, t}$. We identify  $H_0$ with the stalk of $H_{R,\cT}$ at the image of $t$ in $\cT^{\log}$.
Let $\Gamma:=\pi_1(s^{\log})=\pi_1(\cT^{\log})= \Hom(\cS^{\gp}, \Z)$.  Then $\Gamma$ acts on $H_0$, and the local system  $H_{R,\cT}$ has a canonical $\Gamma$-level structure with respect to  the constant sheaf 
$H_0$ (that is, we have a canonical global section  of the quotient sheaf ${\cal I}/\Gamma$ on $\cT^{\log}$, where $\cal I$ is  the sheaf of isomorphisms from  $H_{R,\cT}$ to $H_0$). 
Hence we have a canonical $\Gamma$-level structure on $H_R$ with respect to   $H_0$
on $V^{\log}$ for some open neighborhood $V$ of $s$ in $X$.
We may assume that $X=V$. 
\end{para}

\begin{para}  

The following happens on $\cT^{\log}$ (\cite{KU} 2.3.7).   We can regard $H_0$ as a constant subsheaf of $\cO^{\log}_{\cT} \otimes_R H_{R,\cT}$ as follows. 

Since $\Gamma= \pi_1(\cT^{\log})$ is commutative, $\Gamma$ acts on $H_{R,\cT}$.  

Define a local system $H_0'$ on $\cT^{\log}$ as follows. Taking $(q_j)_{1\leq j\leq n}$ which is a $\Z$-base  of $\cS^{\gp}$,   let $$H'_0:= \xi H_{R,\cT}\subset \cO^{\log}_{\cT}\otimes_R H_{R,\cT}\quad \text{with}\; \xi= \exp\bigl(\sum_{j=1}^n z_jN_j\bigr).$$ 
Here $\xi$ depends on the choices of the branches $z_j$ of $(2\pi i)^{-1}\log(q_j)$ ($1\leq j\leq n$), but   $H'_0$ is independent of the choice.  Furthermore, $\xi \bmod \Gamma$ is independent of the choice of the $\Z$-base $(q_j)_j$ of $\cS^{\gp}$.

Then by \ref{xi1},  $H'_0$ descends to a local system on $\cT$. Since $\cT$ is contractible, $H_0'$ is a constant sheaf. We have an isomorphism 
$$H_0' \overset{\cong}\to H_0$$
by using a ring homomorphism $\cO_{X,t}^{\log}\to \C$ which extends the evaluation $\cO_{X,s}\to \C$ by $z_j\mapsto 0$.

We identify $H'_0$ and $H_0$ via this isomorphism.

We regard $H_0$ as a constant sheaf on  $\cT$ via the above identification. We regard $H_0$ also as a constant sheaf on $X$. We have   $$\cO^{\log}_{\cT} \otimes_{\Z} H_{R,\cT} =  \cO^{\log}_{\cT} \otimes_R H_0, \quad \cO^{\log}_X \otimes_R H_R= \cO^{\log}_X \otimes_R H_0,$$  
$$\tau_*(\cO^{\log}_{\cT} \otimes_{\Z} H_{R,\cT})= \cO_{\cT} \otimes_R H_0,\quad \tau_*(\cO^{\log}_X \otimes_R H_R)= \cO_X \otimes_R H_0. $$
Thus on $X$,  $H_{\cO}=\tau_*(\cO^{\log}_X \otimes_R H_R)$ is identified with $\cO_X \otimes_R H_0$.

Note that in the formula (1) in  \cite{KU} 2.3.7, ``with $\nu=$'' should be replaced by ``with $\xi=$''. 
\end{para}

\begin{para}
\label{Dspcheck}
  Let $h_w^p$ be the $\C$-dimension of the $\gr_F^p$ of $\gr^W_w$ of $H$ at $s$. Note that we have $W$ and $\langle\cdot,\cdot\rangle_w$ on $\gr^W_wH_{0,\Q}$.

Let $\Dc$ be the space of all descending filtrations $F$ on $H_{0,\C}$ such that the rank of $\gr^p_F$ of $\gr^W_w$ is the given $h_w^p$ and such that the annihilator of $F^p\gr^W_{w, \C}$ in 
$\gr^W_{w,\C}$ under $\langle\cdot,\cdot\rangle_w$ is $F^{w+1-p}\gr^W_{w,\C}$.

Then $\Dc$ is a complex analytic manifold. 

The Hodge filtration on $H_{\cO}=\cO_X \otimes_R H_0$ on $X$ gives a morphism $X\to \Dc$, and the Hodge filtration on $H_{\cO}$ is the pullback of the universal Hodge filtration on $\cO_{\Dc} \otimes_R H_0$. Thus we have a morphism $X\to \check{E}:= \cT\times \Dc$. 
  
\end{para}

\begin{para}\label{E1}

 On $\check{E}$, we have the local system, the pullback $H_{R, \check{E}}$ of $H_{R, \cT}$ with $W$ and $\langle\cdot,\cdot\rangle_w$, and we have the Hodge filtration on $\cO_{\check{E}}\otimes_R H_0$ which is the pullback of the universal Hodge filtration of $\cO_{\Dc} \otimes_R H_0$. We have also an isomorphism $\cO_{\check{E}}^{\log}\otimes_R H_{R, \check{E}} \cong \cO_{\check{E}}^{\log} \otimes_R H_0$.
 We denote this object by $H_{\check{E}}$. $H$ on $X$ is the pullback of this $H_{\check{E}}$ under the canonical morphism  (period map) $X\to \check{E}$. 
 \end{para}

\begin{para}

Let $\tilde E$ (resp.\ $E$) be the set of all points $z$ of $\check{E}$ such that the pullback of $H_{\check{E}}$ to $z$ satisfies the Griffiths transversality (resp.\ is an $R$-LMH with polarized $\gr^W$), and  endow  $\tilde E$ (resp.\ $E$) with the strong topology in $\check{E}$ in the sense of \cite{KU} Section 3.1 and with the inverse images of  $\cO_{\check E}$  and the log structure of  $\check E$. 
 Then $E$ is an open set of $\tilde E$, and $E$ and $\tilde E$  are log manifolds.  This is seen by the arguments in \cite{KU} Section 7, \cite{KNU4} Appendix A.1, and \cite{KNU5} 4.5. 

 Let $H_{\tilde E}$ (resp.\ $H_E$) be the pullback of $H_{\check{E}}$ to $\tilde E$ (resp.\ $E$). 
\end{para}

\begin{para}
Since $X\to \check{E}$ is strict, for $x\in X$ with the image $z$ in $\check{E}$, the pullback of $H$ to the fs log point $x$ satisfies the Griffiths transversality (resp.\ is an $R$-LMH with polarized $\gr^W$) if and only if the pullback of $H_{\check{E}}$  to $z$ has the same property.
\end{para}

\begin{para} We prove Proposition \ref{open}. 
  Let $X \to \check{E}$ be as above. 
  Assume that $H$ satisfies the small Griffiths transversality. 
  Then the morphism $X\to \check{E}$  factors through $X\to \tilde E$, and $U$ is the inverse image of $E$. Since $E$ is open in $\tilde E$, $U$ is open in $X$. 
\end{para}

\begin{para}\label{E} We prove Proposition \ref{fromLM}. Let $X \to \check{E}$ be as above. Assume that $H$ is an $R$-LMH with polarized $\gr^W$. Then the morphism  $X\to \check{E}$ factors through $E\subset \check{E}$ and $H$ is  the pullback of $H_E$.
\end{para}

\section{Polarized log mixed Hodge structure and PLH}
\label{s:pLMH}

\begin{para}

Let $X$ be an object of $\cB(\log)$, let $S$ be the standard log point with a fixed generator $q \in M_S$,  let $w\in \Z$, and let $R$ be a subfield of $\R$. Let $\cC_1$ and $\cC_2$ be the following categories. 

Let $\cC_1$ be the category of  pre-$R$-PLH (resp.\  pre-$R$-LVPH) $P$ on $X\times S$ of weight $w$ satisfying the following condition (a).

(a)  Locally on $X$, for some morphism $X\times S \to X \times S$ over $X$, the pullback of $P$ is an $R$-PLH  on the left $X\times S$.

Let $\cC_2$ be the category of $R$-LMH $H$ (resp.\ $R$-LVMH) on $X$ 
endowed with the following structures (i) and (ii) and satisfying the conditions (1) and (2).

(i) A homomorphism $H \otimes H\to R(-w)$ in the category of $R$-LMH such that the induced pairing $\langle\cdot,\cdot\rangle: H_R\times H_R\to R\cdot (2\pi i)^{-w}$ is non-degenerate and $(-1)^w$-symmetric. 

(ii) A homomorphism $N:H\to H(-1)$ such that $N^n: H\to H(-n)$ is zero for some $n\geq 1$ and such that $\langle Nu, v\rangle+\langle u,Nv\rangle =0$.

(1) The weight filtration on $H_R$ coincides with $W(N)[-w]$, where $W(N)$ is the monodromy filtration of $N$.

(2) Let $k\geq w$ and let $\text{Prim}_k$ be the primitive part of $\gr_wH_R$ for $N$. Then the pairing $\text{Prim}_k \times \text{Prim}_k\to R(-k)\;;\; (u,v)\mapsto \langle u, N^{k-w}v\rangle$ is a polarization of the pure PLH $\text{Prim}_k$ of weight $k$.

Morphisms in $\cC_1$ and $\cC_2$ are defined to be isomorphisms in the evident sense.

Note that for an object $H$ of $\cC_2$ and for $k\in \Z$, $\gr^W_kH$ is endowed with  the polarization defined by the decomposition of $\gr^W_kH$ as the direct sum of various primitive parts which are endowed with the polarizations in the condition (2).
\end{para}

\begin{thm}\label{Pmixed0}
We have an equivalence $\cC_1\simeq \cC_2$.

$\cC_1\to \cC_2\;;\; P \mapsto H$
is as follows. Let $\beta$ be the canonical map $(X\times S)^{\log}=X^{\log} \times S^{\log}\to X^{\log}$. 
Then 
$H_R=\beta_*(\exp(\log(q)N)P_R)$, where $N=(2\pi i)^{-1}\log(\gamma)$ for the action $\gamma$ of the canonical generator of $\pi_1(S^{\log})$. 
 $H_{\cO}=P_{\cO}$ with the same Hodge filtration. The isomorphism $\cO_X^{\log}\otimes_R H_R \cong \cO_X^{\log}\otimes_{\cO_X} H_\cO$ is induced from the corresponding isomorphism for $P$. 

$\cC_2\to \cC_1\;; \; H \mapsto P$ 
is as follows.  $P_R= H_R^{(N)}:=\exp(-\log(q)N)H_R\subset \cO_{X\times S}^{\log} \otimes_R H_R$. 
  $P_{\cO}=H_{\cO}$ with the same Hodge filtration. The isomorphism $\cO^{\log}_{X\times S}\otimes_R P_R\cong \cO^{\log}_{X\times S}\otimes_{\cO_X} P_{\cO}$ is induced from the corresponding isomorphism for $H$.

\end{thm}

\begin{para} Here, in the argument of $P\mapsto H$,  the inverse image of $H_R$ on $(X\times S)^{\log}$ is $\exp(\log(q)N)P_R$. In fact, $\gamma$ does not change $\exp(\log(q)N)a$ for an element $a$ of the stalk of $P_R$ because 
$\gamma=(\gamma^*)^{-1}$ ($\gamma^*$ is the pullback by $\gamma$) sends $\log(q)$ to $\log(q)-2\pi i$ and hence $\gamma(\exp(\log(q)N)a)= \exp((\log(q)-2\pi i)N) \exp(2\pi iN)a= \exp(\log(q)N)a$. 
\end{para}

\begin{para}\label{first} We first prove Theorem \ref{Pmixed0} in the case where $X$ is an fs log point and $R=\bR$. 
 In this case, Theorem \ref{Pmixed0} is equivalent to the following Proposition \ref{Pmixed}.  The case $n=0$ of Proposition \ref{Pmixed} is the well-known relation 
between  nilpotent orbits and polarized mixed Hodge structures in \cite{Sc} Theorem (6.16) and in \cite{CKS1} (3.13).

\end{para}

\begin{prop}\label{Pmixed} Let $V$ be a finite dimensional $\R$-vector space. 
Let $w\in \Z$ and let $\langle\cdot, \cdot\rangle:V \times V\to \R$ be a non-degenerate $(-1)^w$-symmetric $\bR$-bilinear form. 
Let $N_0, N_1, \dots, N_n:V\to V$ be  mutually commuting nilpotent linear operators 
such that  $\langle N_ju, v\rangle+\langle u, N_jv\rangle=0$ for all $u,v$ in $V$ and $0\leq j\leq n$.  Let $W=W(N_0)[-w]$, where $W(N_0)$ is the monodromy filtration of $N_0$.  
Let $F$ be a descending filtration on $V_\C$ such that the annihilator of $F^p$ for $\langle\cdot,\cdot\rangle$ is $F^{w+1-p}$ for every $p$.

Then,  $(V, \langle\cdot,\cdot \rangle, N_0, aN_0+N_1, \dots, aN_0+N_n, F)$ generates a pure nilpotent orbit of weight $w$ for any $a\gg 0$  if and only if the following two conditions {\rm (i)} and {\rm (ii)} are satisfied.

{\rm (i)}  $(V, W, N_1, \dots, N_n, F)$ generates a mixed nilpotent orbit.

{\rm (ii)} Let $k \geq w$, let $P_k\subset \gr^W_k$ be the primitive part for $N_0$, and let $\langle\cdot,\cdot \rangle_k: P_k\times P_k\to \R$ be the bilinear form $(u,v)\mapsto \langle u, N_0^{k-w}v\rangle$. 
Then  $(P_k, \langle \cdot, \cdot\rangle_k, N_1, \dots, N_n, F(\gr^W_k)\vert_{P_k})$  is a pure nilpotent orbit of weight $k$.

\end{prop}

\begin{para} The relation between Theorem \ref{Pmixed0} and Proposition \ref{Pmixed} is that $X$ in Theorem \ref{Pmixed0}  provides $N_1, \dots, N_n$ of Proposition \ref{Pmixed} and $S$ in Theorem \ref{Pmixed0} provides $N_0$ of Proposition \ref{Pmixed}.
\end{para}

\begin{para}\label{ifpart} We prove the if part of Proposition \ref{Pmixed}.

By \cite{KNU08} 10.2, we have an action $\tau=(\tau_j)_{0\leq j\leq n}$ of ${\Bbb G}_m^{\{0,\dots, n\}}$ on $V$ associated with the mixed nilpotent orbit  $H=(V, W, N_1, \dots, N_n, F)$  with polarized $\gr^W$ such that  $\tau_0$ splits $W$ and $\tau_j$ for each $1\leq j\leq n$ splits the relative monodromy filtration of  $N_1+\dots+N_j$ with respect to $W$. 
 For $y=(y_0, \dots, y_n)$, $y_j\in \R_{>0}$, let $t(y)= \prod_{j=0}^n \tau_j((y_{j+1}/y_j)^{1/2})$, where $y_{n+1}$ denotes $1$. 

Fix $b:\{0,\dots, n\}\times \{0, \dots, n\}\to [0,\infty]$ such that $b_{j,k}b_{k,\ell}=b_{j,\ell}$ unless the set $\{b_{j,k}, b_{k,\ell}\}$ coincides with the set $\{0,\infty\}$, such that $b_{j,j}=1$ for $0 \leq j \leq n$, and such that $b_{0,j}=\infty$ and $b_{j,0}=0$ for $1\leq j\leq n$.   Then by the $\SL(2)$-orbit theorem for mixed nilpotent orbit (\cite{KNU08}), we have the associated $\hat F$ and $\hat N_1, \dots, \hat N_n$ such that when $y_j\to \infty$ ($0\leq j\leq n$) and  $(y_j/y_k)_{j,k}$ converges to $b$, 
 then $t(y)^{-1}\exp(\sum_{j=0}^n iy_jN_j)F\in \check D$ converges to $\exp(\sum_{j=0}^n i\hat N_j)\hat F\in D$, where $D$ is the classifying space of PH and $\check D$ is its compact dual (see \cite{KNU08} 0.1 for the precise definitions). 
Since $D$ is open in $\check D$, we have $t(y)^{-1}\exp(\sum_{j=0}^n iy_jN_j)F\in D$ and hence $\exp(\sum_{j=0}^n iy_jN_j)F\in D$ if $(y_j/y_k)_{j,k}$ is sufficiently near to $b$. 

This proves that $(V, \langle\cdot,\cdot \rangle, N_0, aN_0+N_1, \dots, aN_0+N_n, F)$ generates a pure nilpotent orbit of weight $w$ for any $a\gg 0$. 

\end{para}

\begin{para}\label{onlyif}  We prove the only if part of Proposition \ref{Pmixed}. There is $c\in \R$ such that if $y_j\geq c$ for $0\leq j\leq n$, $(V, \langle\cdot,\cdot\rangle, F_y)$ with $F_y:=\exp(iy_0N_0+\sum_{j=1}^n iy_j(aN_0+N_j))F=\exp(i(y_0+a\sum_{j=1}^ny_j)N_0+\sum_{j=1}^n i y_jN_j)F$  is a polarized Hodge structure of weight $w$. 
Hence if $b_j\geq c$  for $1\leq j\leq n$, $(V, \langle\cdot,\cdot\rangle, N_0, \exp(\sum_{j=1}^n ib_jN_j)F)$  generates a nilpotent orbit of weight $w$. Hence by the
 classical result of Schmid \cite{Sc} Theorem (6.16) to which we referred to in \ref{first},  if $b_j\geq c$ for $1\leq j\leq n$, $(V, W, \exp(\sum_{j=1}^n ib_jN_j)F)$ is a mixed Hodge structure and for $k\geq w$, $(P_k ,\langle\cdot,\cdot\rangle_k, \exp(\sum_{j=1}^n ib_jN_j)F(\gr^W_k)|_{P_k})$  is a polarized Hodge structure. This proves that  the conditions (i) and (ii) are satisfied. 
\end{para}

\begin{para}  Theorem \ref{Pmixed0} for general $X$ is reduced to the case where $X$ is an fs log point by Proposition \ref{open}.

\end{para}

\section{Study of  Ext groups}
\label{sec:ext}

\begin{para} 
  Let $X$ be an object of $\cB(\log)$. We will consider the following six categories $\cH\supset \cH(*)\supset \cH(**)$, $\cL_0$, $\cL \supset \cL(*)$. 
Fix a subring $R$ of $\R$.

Let $\cH$ (resp.\ $\cH(*)$, resp.\ $\cH(**)$)  be the category of pre-$R$-LMH (resp.\ $R$-LMH, resp.\ $R$-LVMH) on $X$ satisfying the condition (i) in \ref{MP0}.

Let 
 $\cL_0$ be the category of locally constant sheaves  of finite dimensional $\R$-vector spaces on $X^{\log}$ whose local monodromies are unipotent.
  Let $\cL$ be the category of pairs $(L, W)$, where $L$ is an object of $\cL_0$ and $W$ is an increasing filtration (called weight filtration) on $L$ such that 
  each filter $W_k$ is locally constant and that 
$W_k=L$ for some $k$ and $W_k=0$ for some $k$. Let  $\cL(*)$ be the full subcategory of $\cL$ consisting of objects $L$  such that 
 the local monodromies of $L$ are admissible (see \cite{KNU3} 1.2.4).

These categories are exact categories (in the sense of Quillen). For a contemporary treatment of  exact categories, see \cite{B}.  A short exact sequence in $\cH$, in $\cH(*)$, or in $\cH(**)$  is a sequence $0\to H_1\to H_2\to H_3\to 0$ such that 
$0\to W_kH_{1,R}\to W_kH_{2,R}\to W_kH_{3,R}\to 0$ for all $k$ and  $0\to F^pH_{1,\cO}\to F^pH_{2,\cO}\to F^pH_{3,\cO}\to 0$ for all $p$ are exact. A short exact sequence in $\cL_0$ is an evident one, and that in $\cL$ or in $\cL(*)$ is a sequence $0\to L_1\to L_2\to L_3\to 0$ such that the sequences $0\to W_kL_1\to W_kL_2\to W_kL_3\to 0$ are exact for all $k$.

We have Yoneda's  higher Ext group $\Ext^n$ for any exact category (\cite{P} Appendix Proposition A.13). For example, $\Ext^1$ is 
 the set of isomorphism classes of extensions endowed with the group law given by Baer sums. A short exact sequence gives a long  exact sequence of $\Ext^n$.

\end{para}

\begin{para} If $\cC$ is one of the above categories 
$\cH$, $\cH(*)$, $\cH(**)$, $\cL_0$, $\cL$, and $\cL(*)$, 
and if $A, B$ are objects of $\cC$, we have a sheaf $\cE xt^n_{\cC}(A, B)$ of abelian groups on $X$ which is the sheafification of the presheaf $U\mapsto \Ext^n_{\cC'}(A',B')$, where  $U$ is an open set of $X$ and $\cC', A', B'$ are the restrictions of $\cC$, $A$, $B$ over $U$, respectively.

\end{para}

The goal of this Section \ref{sec:ext} is to prove the following proposition.

\begin{prop}\label{propext2} Let $\cC$ be either $\cH(*)$ or $\cH(**)$. 
 Let 
$0\to R(1) \to P \to Q \to 0$ be an exact sequence in $\cC$ and assume that  the weights of $Q$ $\leq -1$. Then the map $\cExt^1_{\cC}(R, P) \to \cExt^1_{\cC}(R, Q)$ is surjective.

\end{prop}

\begin{lem}\label{L(adm)} 
Let  $L$ be an object of $\cL(*)$ such that the local monodromy actions on $L$ are trivial and such that $W_{-2n}L=0$.

Then 
the canonical map $\cE xt^n_{\cL(*)}(\R, L) \to \cE xt^n_{\cL_0}(\R, L) = R^n\tau_*L$ 
is the zero map. Here $W_0\R=\R$ and $W_{-1}\R=0$. 
 
In other words (assume $n\geq 1$), if  $$0\to L\to L_n\to \dots \to L_2 \to L_1\to  \R \to 0$$ is an exact sequence 
in $\cL(*)$,  the induced section of $R^n\tau_*L$ is $0$.

\end{lem}

\begin{pf} The case $n=0$ is evident. Assume $n\geq 1$ and consider an exact sequence $0\to L\to L_n\to \dots \to L_2 \to L_1\to  \R \to 0$ in $\cL(*)$. 
Let $t\in X^{\log}$ and let $x\in X$ be the image of $t$. 
Then  $(R^n\tau_*L)_x$ is
 isomorphic to $(\bigwedge^n_{\Z}\;  (M_X^{\gp}/\cO_X^\times)_x)\otimes_{\Z} L_t$ and is isomorphic to the Lie cohomology $ H^n(\frak g, L_t)$, where $\frak g$ is the commutative Lie algebra $\Hom((M_X^{\gp}/\cO^\times_X)_x,  \Q)$ which acts on $L_t$ trivially. 
Let $\sig$ be the monodromy cone 
$\Hom((M_X/\cO^\times_X)_x, \R_{\geq 0}^{\text{add}})$ of $x$. 
  In the rest, we omit $t$ in $(\cdot)_t$.
  By the admissibility, we have a relative monodromy filtration $W(\sig)$ on $L$, $L_j$ ($1\leq j\leq n$), $\R$ and the sequence $0\to W(\sig)_kL \to W(\sig)_k L_n\to \dots \to W(\sig)_kL_1 \to W(\sig)_k\R \to 0$ is exact for every $k$. 
For $0\leq j\leq n$, let $I_j$ be the image of $L_{j+1} \to L_j$, where $L_0$ denotes $\R$ and $L_{n+1}$ denotes $L$. Hence $I_0$ is identified with $\R$ and $I_n$ is identified with $L$.  We compute the $n$-cocycle $f_n: \bigwedge^n_{\Q} \frak g\to L$ corresponding to the element of $H^n(\frak g, L)$ in problem, by the standard method: By induction on $j$, we get a $j$-cocycle $\bigwedge^j_\Q \frak g \to I_j$ starting from $f_0 =1\in \R=I_0$ . To get $f_{j+1}$ from $f_j$, we lift $f_j$ to a $j$-cochain $\tilde f_j: \bigwedge^j_\Q \frak g\to L_{j+1}$ and obtain $f_{j+1}$ from $\tilde f_j$. By induction on $j$, we get $f_j$ whose image is in $W(\sig)_{-2j}I_j$ and lift it to $\tilde f_j$ whose image is in $W(\sig)_{-2j}L_{j+1}$, and get $f_{j+1}$ whose image is contained in $(\sum_{N\in \frak g} NW(\sig)_{-2j}L_{j+1})\cap I_{j+1}\subset  W(\sig)_{-2(j+1)}I_{j+1}$. Thus we get $f_n$ whose image is in $W(\sig)_{-2n}L=W_{-2n}L=0$. 
\end{pf}

\begin{lem}\label{adm} Let $0\to L \to P\to Q\to 0$ be an exact sequence in $\cL(*)$. Assume that the local monodromy actions on $L$ are trivial and  $W_0L=L$. For $a\in \cE xt^1_{\cL}(\R, P)$ and the image $b$ of $a$ in $\cE xt^1_{\cL}(\R, Q)$, $a$  belongs to $\cE xt^1_{\cL(*)}(\R, P)$ if and only if $b$ belongs to $\cE xt^1_{\cL(*)}(\R, Q)$. 

\end{lem}

\begin{pf}  It is enough to prove the if part. We may assume that $a$ is the class of an exact sequence $0 \to P \to \tilde P \to \R \to 0$ in $\cL$ and $b$ is the class of the exact sequence  $0\to Q \to \tilde Q \to \R \to 0$ in $\cL(*)$, where $\tilde Q=\tilde P/L$.
  Let $\sig$ be as in the proof of Lemma \ref{L(adm)}. 
For each face $\sig'$ of $\sig$, we have the relative monodromy filtration $W(\sig')$ on $\tilde P$ defined as follows. If $k\leq -1$, $W(\sig')_k\tilde P= W(\sig')_kP$. If $k\geq 0$, $W(\sig')_k\tilde P$ is the inverse image of $W(\sig')_k\tilde Q$. 
  Hence $\tilde P$ belongs to $\cL(*)$. 
\end{pf}

\begin{para}\label{top} Note that for a topological space $T$ and for a complex of sheaves 
of abelian groups on $T$ of the form $C=[C^0\to C^1]$ (that is, a complex of sheaves of abelian groups whose degree $d$-parts are zero unless $d=0, 1$), the hyper-cohomology $H^1(T, C)$ is identified with the set of isomorphism classes of pairs 
of an exact sequence of the form $0\to C^0\to E \to \Z\to 0$ and a splitting $E\to C^1$.
\end{para}

\begin{para}\label{C(H)}

For a pre-$R$-LMH $H$ satisfying the condition (i) in \ref{MP0} 
on $X$, let 
$$C(H)=[H_R\to (H_R\otimes \cO_X^{\log})/F^0],$$
where $H_{R}$ is of degree 0. 
  Then we have identifications
$$\Hom_{\text{pre-$R$-LMH}/X}(R, H) = H^0(X^{\log}, C(H)) \quad\text{if the weights of $H$ $\leq 0$},$$
$$ \Ext^1_{\text{pre-$R$-LMH}/X}(R, H) = H^1(X^{\log}, C(H))\quad\text{if  the weights of $H$ $\leq -1$}.$$
This is shown by using \ref{top}. 
\end{para}

\begin{lem}\label{extlem1} Let $0\to R(1) \to P \to Q\to 0$ be an exact sequence in $\cH$ such that the weights of $Q$ $\leq -1$. Then we have an exact sequence 
$$\cE xt^1_{\cH}(R, P) \to \cE xt^1_{\cH}(R, Q) \to R^2\tau_*R(1).$$

\end{lem}

\begin{pf} By $\cE xt^1_{\cH}(R, P)\cong R^1\tau_*C(P)$ (\ref{C(H)}) and by the corresponding isomorphism for $Q$, we have an exact sequence
$\cE xt^1_{\cH}(R, P)\to \cE xt^1_{\cH}(R, Q) \to R^2\tau_*C(R(1))$. Since $C(R(1))=[R(1) \to \cO_X^{\log}]$, we have an exact sequence 
$0= R^1\tau_*\cO^{\log}_X \to R^2\tau_*C(R(1)) \to R^2\tau_*R(1)$. 
\end{pf}

\begin{lem}\label{extlem2}  Let $\cC$ and $0 \to R(1) \to P \to Q \to 0$ be as in the hypothesis of Proposition $\ref{propext2}$. 
  Then the map $\cE xt^1_{\cH}(R, Q) \to R^2\tau_*R(1)$ in Lemma $\ref{extlem1}$ kills $\cE xt^1_{\cC}(R, Q)$. 

\end{lem}

\begin{pf} Since $R^2\tau_*R(1)\cong R\otimes_{\Z} \bigwedge^2_{\Z}\;  (M^{\gp}_X/\cO_X^\times)$, the map  $R^2\tau_*R(1)\to R^2\tau_*\R(1)$ is injective.
The map $\cE xt^1_{\cC}(R,Q)\to R^2\tau_*\R(1)$ factors as 
$\cE xt^1_{\cC}(R, Q)\to \cE xt^1_{\cL(*)}(\R, Q_\R)\to \cE xt^2_{\cL(*)}(\R, \R(1)) \to \cE xt^2_{\cL_0}(\R, \R(1))=R^2\tau_*\R(1)$. This  is $0$ by the case $n=2$ and $L=\R(1)$ of Lemma \ref{L(adm)}.
\end{pf}

 By Lemmas \ref{extlem1} and \ref{extlem2}, for the proof of Proposition \ref{propext2}, it is sufficient to prove the following. 

\begin{lem}\label{extlem3} Let  $\cC$ and $0 \to R(1) \to P \to Q \to 0$  be as in the hypothesis of Proposition $\ref{propext2}$. 
  Then for $a\in \cE xt^1_{\cH}(R, P)$ and the image $b$ of $a$ in $\cE xt_{\cH}^1(R, Q)$, $a$ belongs to $\cE xt^1_{\cC}(R, P)$ if and only if $b$ belongs to $\cE xt^1_{\cC}(R, Q)$.

\end{lem}

\begin{pf} It is sufficient to prove the if part. We may assume that $a$ is the class of an exact sequence $0 \to P \to \tilde P \to R \to 0$ in $\cH$ and $b$ is the class of the exact sequence $0\to Q \to \tilde Q  \to R \to 0$ in $\cH(*)$, where $\tilde Q=\tilde P/R(1)$.

Then $\tilde P$ satisfies the admissibility of local monodromy by  Lemma \ref{adm}.

If $\cC=\cH(*)$ (resp.\ $\cC=\cH(**)$), $\tilde P$ satisfies the small (resp.\ big) Griffiths transversality because
$$((\tilde P_R\otimes_R \cO^{\log}_x)/F^{-1}) \otimes_{\C} \omega^1_x\overset{\cong}\to ((\tilde Q_R\otimes_R \cO^{\log}_x)/F^{-1}) \otimes_{\bC} \omega^1_x\quad \text{for $x\in X$}.$$

$${\text (\rm{resp.}\ } ((\tilde P_R\otimes_R \cO^{\log}_X)/F^{-1}) \otimes_{\cO_X} \omega^1_X\overset{\cong}\to ((\tilde Q_R\otimes_R \cO^{\log}_X)/F^{-1}) \otimes_{\cO_X} \omega^1_X.)$$
  Finally, we prove that for $t\in X^{\log}$ lying over $x\in X$ and for 
$s\in \spe(t)$, if $\exp(s(\log(\cdot))):M_{X,x}^{\gp} \to \bC^{\times}$ is  
sufficiently near to the structure homomorphism $\alpha_{X,x}$, then the associated specialization $\tilde P(s)$ is a mixed Hodge structure. 
If $k\neq 0$, $\gr^W_k \tilde P=\gr^W_kP$  and hence $\gr^W_k\tilde P(s)$ is a  Hodge structure of weight $k$.
 We consider $\gr^W_0\tilde P(s)$. We have exact sequences $0\to \gr^W_0P_\C \to \gr^W_0\tilde P_\C \to \C\to 0$, $0\to F^p\gr^W_0P(s)_\C \to F^p\gr^W_0\tilde P(s)_\C\to \C\to 0$ for $p\leq 0$ and $F^p\gr^W_0P(s)_\C\overset{\cong}\to F^p\gr^W_0\tilde P(s)_\C$ for $p\geq 1$. 
Hence we have the Hodge decomposition of  $\gr^W_0\tilde P(s)$. 
\end{pf}

\section{Construction of an extension with three graded quotients}\label{sec:three}

  In this section, we prove Theorem \ref{MP} in the special case where $W_{-1}H=H$,  $W_{-3}H=0$, and $\gr^W_{-2}H=R(1)$. Let $Q:=\gr^W_{-1}H$. Thus we have an exact sequence $0\to R(1) \to H \to Q\to 0$. 

\begin{para} Let $\tilde Q= H^*(1)$, where $(\cdot)^*$ denotes the dual. The intersection form $\langle\cdot, \cdot \rangle_{-1}: Q \times Q \to R(1)$  of a polarization gives an isomorphism $Q\overset{\cong}\to  Q^*(1)\;;\;u\mapsto (v\mapsto \langle u,v\rangle)$. Hence $W_0\tilde Q=\tilde Q$, $\gr^W_0\tilde Q= R$, $\gr^W_{-1}\tilde Q=Q$, and $W_{-2}\tilde Q=0$.

Let $\cC$ be $\cH(*)$ or $\cH(**)$, and assume that $H$ belongs to $\cC$. By the surjectivity of $\cE xt^1_{\cC}(R, H)\to \cE xt^1_{\cC}(R, Q)$ (Proposition \ref{propext2}), locally on $X$, the class of $\tilde Q$ in $\Ext^1_{\cC}(R, Q)$ lifts to the class in $\Ext^1_{\cC}(R,H)$ of an exact sequence $0\to H \to \tilde H \to R \to 0$ such that $\tilde H/R(1)=\tilde Q$. 

\end{para}

\begin{lem}\label{dual} There is a unique isomorphism $\tilde H\cong (\tilde H)^*(1)$ whose $\gr^W_{-1}$ coincides with $Q\cong Q^*(1)$, whose $\gr^W_0$ is the identity homomorphism of $R$, and  whose $\gr^W_{-2}$ is the multiplication by $-1$ on $R(1)$. 
\end{lem}

\begin{pf} 
If we have two such isomorphisms $f,g:\tilde H\cong (\tilde H)^*(1)$, $f-g: \tilde H \to (\tilde H)^*(1)$ sends $W_k$ to $W_{k-1}$ and hence is zero. This proves the uniqueness. 

We prove the existence. 

First we consider the case $Q=0$. We have an exact sequence $0\to R(1)\to \tilde H_R \to R\to 0$. We have the wedge product  $\langle\cdot,\cdot\rangle: \tilde H_R \times \tilde H_R \to \bigwedge^2_R \tilde H_R\cong R(1)$, where the last isomorphism is  such that the induced map
$\gr^W_0\otimes \gr^W_{-2} \to R(1)$ is the canonical map $R \otimes R(1)\to R(1)$. This pairing induces a perfect pairing $\tilde H_{\cO}\times \tilde H_{\cO}\to \cO_X$ of $\cO_X$-modules. The Hodge filtration $F$ of $H_{\cO}$ satisfies  $F^{-1}=H_{\cO}$, $F^0$ is a line bundle, $F^1=0$,  and the annihilator of $F^0$ under $\langle\cdot,\cdot\rangle$ is $F^0$. Hence $\langle\cdot,\cdot\rangle$ induces  $\tilde H \cong (\tilde H)^*(1)\;;\; u\mapsto (v\mapsto \langle u, v\rangle)$.

We consider the general case.
A proof may be given by  writing explicitly the involved Hodge filtrations, but we give a proof by some abstract argument. 
We may assume $\cC=\cH(*)$.

 We have two sections $a,a^*(1)$ of $\cE xt^1_{\cH(*)}(R, H)$ given by $\tilde H$ and $(\tilde H)^*(1)$, respectively. Here
to define $a^*(1)$, we  identify
$W_{-1}((\tilde H)^*(1))$ with $H$ via 
$W_{-1}((\tilde H)^*(1)) = (\tilde Q)^*(1) = (H^*(1))^*(1)=H$.

For the proof of Lemma \ref{dual}, it is sufficient to prove

\medskip

(1) $a+a^*(1)=0$. 

\medskip

 In fact, then, there is a natural isomorphism from $\tilde H$ to the pushout of $(\tilde H)^*(1)$ by $-1: W_{-1}((\tilde H)^*(1)) \to W_{-1}((\tilde H)^*(1))$, whose $\gr^W_0$ is the identity homomorphism of $R$ such that the induced $H=W_{-1}(\tilde H) \to W_{-1}((\tilde H)^*(1)) =H$ is the multiplication by $-1$.
  This is the desired isomorphism because $Q \cong Q^*(1) \cong (Q^*(1))^*(1)=Q$ is the  multiplication by $-1$ on $Q$.

  Note that we have already proved (1) in the case $Q=0$, by the above consideration on the case $Q=0$. 

Concerning (1), we have the following (2) and (3).

\medskip

(2) This element $a+a^*(1)$  belongs to 
$\cE xt^1_{\cH(*)}(R, R(1))$, since its image in $\cE xt^1_{\cH(*)}(R,Q)$ is zero. 

\medskip

(3) 
This sum $a+a^*(1)$ is determined by the image of $a$ in $\cE xt^1_{\cH(*)}(R, Q)$. In fact, if the image of $b\in \cE xt^1_{\cH(*)}(R, H)$ in $\cE xt^1_{\cH(*)}(R, Q)$ coincides with that of $a$, we have $b=a+u$ for some  $u\in \cE xt^1_{\cH(*)}(R, R(1))$. By the case $Q=0$ considered above, $u+u^*(1)=0$, and hence $b+b^*(1)=(a+u)+(a+u)^*(1)=(a+a^*(1))+ (u+u^*(1))= a+a^*(1)$.  

\medskip

By (2) and (3), 
 we obtain a homomorphism 
$\cE xt^1_{\cH(*)}(R, Q)\to \cE xt^1_{\cH(*)}(R, R(1))$. Note that $\cE xt^1_{\cH(*)}(R, R(1))\subset \cE xt^1_{\cH}(R, R(1))= \tau_*(\cO_X^{\log}/R(1))$. 
We prove that this homomorphism is the zero map. By Proposition \ref{fromLM}, this is reduced to the case where $X$ is  a log manifold. 
Let $X_{\triv}$ be the open set of $X$ consisting of all points at which the log structure of $X$ is trivial, which is a  complex analytic manifold. In this case, for the inclusion map $j^{\log}: X_{\triv}\to X^{\log}$, $\cO_X^{\log} \to j^{\log}_*(\cO_{X_{\triv}})$  is injective and $R(1) \to j^{\log}_*(R(1))$ is an isomorphism, and hence the map $\cO_X^{\log}/R(1) \to j^{\log}_*(\cO_{X_{\triv}}/R(1))$ is injective. By this, replacing $X$ by $X_{\triv}$, we are reduced to the case where $X$ is 
a complex analytic manifold with the trivial log structure. In this case, $\cH(*)=\cH$. 

We first consider the case $X=\Spec(\C)$ with the trivial log structure and $R=\R$. In this case,   
$\Ext^1_{\cH}(\R, Q)=Q_\R\bs Q_{\C}/F^0=0$ because $Q$ is of weight $-1$, and hence the map $\Ext^1_{\cH}(\R, Q)\to \Ext^1_{\cH}(\R, \R(1))$ is the zero map. 

This proves that in the case where $X$ is a complex analytic manifold with the trivial log structure, the composite map $\cE xt^1_{\cH}(R, Q)\to  \cE xt^1_{\cH}(R, R(1))= \cO_X/R(1) \to \cO_X/\R(1)$ is the zero map. Hence in this case, the map 
$\cE xt^1_{\cH}(R, Q)=Q_R\bs (Q_R \otimes_R \cO_X)/F^0\to \cE xt^1_{\cH}(R, R(1))= \cO_X/R(1)$ has the image in $\R(1)/R(1)$. Since there is no non-zero homomorphism from $\cO_X$ to the constant sheaf $\R(1)/R(1)$ which is functorial in $X$,  we  see that this map is zero. 
\end{pf}

\begin{para} 
  We return to the proof of Theorem \ref{MP}. 
  Let $N: \tilde H\to \tilde H(-1)$ be the homomorphism $\tilde H \to \gr^W_0\tilde H\cong R  \cong (\gr^W_{-2}\tilde H)(-1) \subset \tilde H(-1)$. 
  Then the desired PLH  (resp.\ VPLH)  on $X\times S$ of weight $-1$ is the pullback $H'$ of $\tilde H^{(N)}$ by $X\times S \to X \times S$ over $X$, $q$ on the right hand side is sent to $qf$ on the left hand side, where $f\in M_X$. 
  Here, $(\cdot)^{(N)}$ is as in Theorem \ref{Pmixed0}.  (Although 
$R$ is a field there, the construction is the same.  
  Note that since $2 \pi i N\tilde H_R\subset \tilde H_R$ and $N^2=0$ here, we have $\exp(2\pi i N)\tilde H_R=\tilde H_R$ and hence $\tilde H^{(N)}$ has an $R$-structure.)
  We prove that this $H'$ 
is a PLH  (resp.\ VPLH). 
  The small or big Griffiths transversality is easy to see. We have the isomorphism $H' \cong (H')^*(1)$ in Lemma \ref{dual} and it gives a polarization $\langle\cdot,\cdot\rangle: H' \times H'\to R(1)$ for which the isomorphism is $u\mapsto (v\mapsto \langle u, v\rangle)$.   To prove that $\langle\cdot,\cdot\rangle$ is actually a polarization, we may assume that the base $X$ is an fs log point and hence we reduce to the case of Theorem \ref{Pmixed0} where $X$ is an fs log point and $R=\R$.  
\end{para}

\section{Proof of Theorem \ref{MP}}\label{sec:pf}

We prove the general case of Theorem \ref{MP}. We use the induction on $n$ for the integer $n\geq 1$ such that there is $w\in \Z$ satisfying $W_wH=H$ and $W_{w-n}H=0$.

\begin{para}\label{MPpf1}  First assume $n=2$. Let $A=\gr^W_wH$, $B=\gr^W_{w-1}H=W_{w-1}H$. We consider $B^*(1) \otimes H$ and use $ B^*(1) \otimes B \to R(1)$ to get $P$ as the pushout. We have $W_{-1}P=P$, $W_{-3}P=0$, $\gr^W_{-1}P=B^*(1)\otimes A$, $\gr^W_{-2}P=R(1)$. 
We get $\tilde P$ as in the previous section. Let $\tilde H=B(-1)\otimes \tilde P$. Since $\tilde P$ is regarded as a PLH  (resp.\ VPLH) of weight $-1$ by the previous section (after taking $(\cdot)^{(N)}$), $\tilde H$ is regarded as a PLH  (resp.\ VPLH) of weight $w$. This $\tilde H$ regarded as a PLH  (resp.\ VPLH) is the desired  $H'$.

Here $N$ is the evident isomorphism   $\gr^W_{w+1}\to \gr^W_{w-1}(-1)$ induced by the identity map of $B$.

We embed $H$ into $\tilde H$  by
$$H \to B \otimes B^*\otimes H  \to B(-1) \otimes P \to B(-1) \otimes \tilde P=\tilde H.$$

The map $H_R \to \tilde H_R$  is injective (with locally free cokernel) because its $\gr^W_{w-1}$ is $B_R \to B_R \otimes B_R^*\otimes B_R \to B_R$, where the first arrow is $b \mapsto \sum_j e_j \otimes e_j^*\otimes b$ for a local base $(e_j)_j$ of $B_R$ and for the dual base $e^*_j$ of $B^*_R$ and the second arrow is the map  $u\otimes v\otimes b \mapsto uv(b)$, and the composite map $B_R\to B_R$ is the identity map because it gives $e_k\mapsto \sum_j e_j\otimes e_j^*\otimes e_k\mapsto e_k$.

\end{para}

\begin{para}\label{MPpf2}  Assume $n\geq 3$, $W_wH=H$, $W_{w-n}H=0$.  Then by the hypothesis of induction applied to $W_{w-1}H$, we find a PLH  (resp.\ an LVPH)  $I\supset W_{w-1}H$ on $X\times S$ of weight $w-1$. 
Let $J$ be the pushout of $H \leftarrow W_{w-1}H \to I$. Then we have $W_wJ=J$, $W_{w-2}J=0$, $\gr^W_wJ=\gr^W_wH$, $\gr^W_{w-1}J= I$. Hence we  find a PLH (resp.\ an LVPH) $K\supset J$ on $X\times S\times S$ of weight $w$.  We get the desired PLH (resp.\ LVPH) $H'\supset H$ on $X\times S$ of weight $w$ as the pullback of $K$ by the diagonal $S\to S\times S$.

It is clear that the condition (a) in Theorem \ref{MP} is satisfied. In \ref{MPpf3} below, we  show that the condition (b) in Theorem \ref{MP} is satisfied.

This  will complete the proof of Theorem \ref{MP}. 

\end{para}

\begin{para}\label{MPpf3}  We prove that the condition (b) is satisfied.   We may assume that $X$ is an fs log point. 
We can use the following. 

If $(W, N_1, \dots, N_n, F)$ is a mixed nilpotent orbit with polarizable $\gr^W$, we have:

(1) We have a relative monodromy filtration $W^{(1)}$ of $(W, N_1)$ and $(W^{(1)}, N_2, \dots, N_n, F)$ is a mixed nilpotent orbit.

(2) The relative monodromy filtration $W^{(2)}$ of $(W^{(1)}, N_2)$  coincides with the relative monodromy filtration of $(W, N_1+N_2)$. 

(3) $(W, N_1+N_2, N_3, \dots, N_n,F)$ is a mixed nilpotent orbit. 

In \ref{MPpf2}, we apply these to the pure weight filtration $W$ of weight $w$ and to the mixed nilpotent orbit 
$K=(W, N_{-1}, N_0, N_1, \dots, N_n, F)$, where $N_1, \dots, N_n$ comes from $X$, $N_0$ comes from $S$ of $X\times S$, and $N_{-1}$ comes from the 
second $S$ of $X\times S \times S$. Then by (1) and (2), we have a mixed nilpotent orbit $(W^{(2)}, N_1, \dots, N_n, F)$, where $W^{(2)}$ is the relative monodromy filtration of $(W, N_{-1}+N_0)$. $I$ (resp.\ $J$) is a sub mixed nilpotent orbit which is the weight $w-1$  (resp.\ $\leq w$) part of 
the mixed nilpotent orbit  $(W^{(1)}, N_0, \dots, N_n, F)$ (consider (1)), where $W^{(1)}$ is the relative monodromy filtration 
of $(W, N_{-1})$. 
We have a mixed nilpotent orbit $I'$ (resp.\ $J'$) $=(W^{(2)}, N_1, \dots, N_n, F)$ whose underlying space is the same as $I$  (resp.\ $J$).

 The original $W_{w-1}H$ is a sub mixed nilpotent orbit of $I'$. Hence the weight filtration of $W_{w-1}H$ is the restriction of $W^{(2)}$. The original $H$ satisfies that $H/W_{w-1}H$ is of weight $w$. $J' /I'$ also has weight $w$. Hence  the weight filtration of $H$ is the restriction of $W^{(2)}$. 
  Thus we complete the proof of that the condition (b) is satisfied. 
\end{para}

\appendix\section{Appendix. Mixed motives and log pure motives}
\label{s:motive}

In this section, we will give some simple construction of the category of mixed motives over a field 
based on the idea that mixed motives should be embedded into log pure motives.

 \begin{para}\label{Hodge}  We present the  Hodge analogue of our story on motives.

Here, Hodge or log Hodge structure means $\Q$-Hodge or $\Q$-log Hodge structure. 

Let $S$ be the standard log point over $\C$. Let (LH) be the category of polarizable log Hodge structures on $S$. For $H\in$ (LH), 
 let $H^{\flat}$ be the associated mixed Hodge structure endowed with the monodromy operator $N:H^{\flat}\to H^{\flat}(-1)$. 

Let (LH$\flat$) be the category of pairs of a mixed Hodge structure and $N$ of the form $H^{\flat}$ with $H\in $(LH). 

Let (MH)  be the category of contra-variant functors from (LH$\flat$) to the category of $\Q$-vector spaces defined by a pair $(H, V)$ as in (1) below, where $H\in$ (LH$\flat$)  and $V$ is a $\Q$-subspace of $H_\Q$ satisfying the following conditions (i) and (ii).

(i)  For some $H'\in$ (LH$\flat$) and for some morphism $H'\to H$ in (LH$\flat$), $V$ is the image of $H'_{\Q}\to H_\Q$. 

(ii) $N$ of $H$ kills $V$. 

\medskip

(1) $H'\mapsto \{h\in  \Mor_{(\LH\flat)}(H', H)\;|\; h(H'_\Q)\subset V\}$. 

\medskip

 Then  (MH) is equivalent to the category of mixed Hodge structures with polarizable $\gr^W$ by the case $X=\Spec(\C)$ with the trivial log structure of  Theorem \ref{MP} and Remark \ref{MPrem} (1).
\end{para}

\begin{para} 

Let $k$ be a field, and let $S$ be the standard log point over $k$.

 We have the log absolute Galois group $\pi_1^{\log}(S)$. It is the automorphism group of the log separable closure 
(cf.\ \cite{Nc} (2.5)) 
$\overline S$ of $S$ over $S$. We have an exact sequence
 $$0\to \hat \Z(1)' \to \pi_1^{\log}(S)\to \Gal(\overline k/k)\to 1,$$
 where $\hat \Z(1)' $ is the product of $\Z_\ell(1)$ for all prime numbers $\ell$ which are not equal to the characteristic of $k$ and $\overline k$ is the separable closure of $k$. 
 \end{para}

 \begin{para} We fix a prime number $\ell$ which is different from the characteristic of $k$. Let $\cP$ be the category of projective vertical log smooth saturated fs log schemes over $S$ which have charts of the log structure Zariski locally. 

Let $X\in \cP$. For $m\in \Z$, let $$H^m(X)_{\ell}=H^m_{\text{log\'et}}(X \times_S {\overline S}, \Q_\ell).$$ 
It is a finite dimensional $\Q_\ell$-vector space endowed with a continuous action of $\pi_1^{\log}(S)$. 

We make symbols $H^m(X)(r)$ and $H^m(X)(r)^{\flat}$ for $X\in \cP$ and for $m, r\in \Z$ such that $m\geq 0$. 
 
 \end{para}

 \begin{para}\label{val} Let $K=K_n$ or $K=KH_n$, where $K_n$ is Quillen's K-theory and $KH_n$ is the homotopy $K$-theory of Weibel \cite{W1}. Following the method in \cite{IKNU} 2.4.6, for an fs log  scheme $X$ having charts Zariski locally, we define $$K_{\lim}(X) =\varinjlim K(X'),$$ where $X'$ ranges over all log modifications of $X$ in the sense of \cite{IKNU} 2.3.6 and $K(X')$ means $K$ of the underlying scheme of $X'$.
The $K$-group $K_{0,\lim}$ is used in \cite{IKNU} and also in the first half of this Appendix, and the $K$-group $KH_{n,\lim}$ is used in the latter half of this Appendix. 
 \end{para}
 
  \begin{para}\label{mor}  Let $X,Y\in \cP$.

 By a morphism $H^m(X)(r)^{\flat} \to H^n(Y)(s)^{\flat}$, we mean a $\Q_{\ell}$-linear map
 $H^m(X)_{\ell}(r)\to H^n(Y)_{\ell}(s)$
which is obtained as below from an element of 
$$\gr^u K_{0,\lim}(X\times_S Y \times {\mathbb G}_m^t)\otimes \Q,$$  
where $t=(n-2s)-(m-2r)$ and $u=d+n-m+r-s$ with $d$ being the dimension of $X$ (the dimension is defined as a locally constant function on $X$), and $\gr^u$ is the graded quotient for the $\gamma$-filtration. (If $X$ is not equi-dimensional, this $K$-group is defined as the direct sum of the $K$-group of connected components of $X$ by using the dimension of each connected component.)

 If $m-2r> n-2s$, there is no non-zero morphism. We assume $m-2r\leq n-2s$. 
 We have homomorphisms
 $$\gr^u K_{0,\lim}(X\times_S Y\times {\mathbb G}_m^t)\otimes \Q \to H^{2u}(X\times_S Y\times {\mathbb G}_m^t)_{\ell}(u)$$ $$ \to H^{2d-m}(X)_{\ell} \otimes H^n(Y)_{\ell}(d+s-r)\otimes H^t({\mathbb G}_m^t)_{\ell}(t)
 \to \Hom(H^m(X)_{\ell}(r), H^n(Y)_{\ell}(s)).$$
Here to have the first homomorphism, we use the fact that the log blowing-up along the log structure does not change the log \'etale cohomology. The second homomorphism is by K\"unneth formula, and the third one is by Poincar\'e duality and by the canonical map $H^t({\mathbb G}_m^t)_{\ell}(t)\overset \cong \to \Q_\ell$ 
induced by $H^1({\mathbb G}_m)_{\ell}\cong \Q_\ell(-1)$.
  (For basic properties of log \'etale cohomology, see \cite{Nc2}.)  
  
As is easily seen,  a linear map $H^m(X)_{\ell}(r)\to H^n(Y)_{\ell} (s)$ commutes with the action of $\pi_1^{\log}(S)$ if it is a morphism $H^m(X)(r)^{\flat}\to H^n(Y)(s)^{\flat}$.

By a  morphism $H^m(X)(r) \to H^n(Y)(s)$, we mean  a morphism $H^m(X)(r)^{\flat}\to H^n(Y)(s)^{\flat}$ such that we can take $t=0$ in the above. 

If $m-2r\neq n-2s$, there is no morphism $H^m(X)(r) \to H^n(Y)(s)$.

\end{para}

\begin{prop}\label{comp} $(1)$ The identity map of $H^m(X)_{\ell}(r)$ is a morphism $H^m(X)(r)\to H^m(X)(r)$ and hence a morphism $H^m(X)(r)^{\flat}\to H^m(X)(r)^{\flat}$. 

$(2)$ For morphisms $H^{m(1)}(X_1)(r(1))^{\flat} \to H^{m(2)}(X_2)(r(2))^{\flat}$ and 
$H^{m(2)}(X_2)(r(2))^{\flat} \to H^{m(3)}(X_3)(r(3))^{\flat}$, 
the composition is a morphism $H^{m(1)}(X_1)(r(1))^{\flat}\to H^{m(3)}(X_3)(r(3))^{\flat}$.  The non-$\flat$ version is also true.
\end{prop}

\begin{pf} The proof for the non-$\flat$ version is given in \cite{IKNU} Propositions 3.1.4 and 3.1.6.
  The $\flat$ version is proved in the same way. 
\end{pf}

Thus we have the category (\LM$\flat$) of $H^m(X)(r)^{\flat}$  and the category (\LM) of $H^m(X)(r)$. The latter was considered in \cite{IKNU}.

\begin{para} We define the category (\MM) as the category of contra-variant functors from (\LM$\flat$) to the category of $\Q$-vector spaces which are obtained as in (1) below, from an object $H^n(Y)(s)^{\flat}$ of  (\LM$\flat$) and a  $\Q_{\ell}$-subspace $V$ of $H^n(Y)_{\ell}(s)$ satisfying the following conditions (i) and (ii).

(i) There is a morphism $H^m(X)(r)^{\flat}\to H^n(Y)(s)^{\flat}$ for some $X, m, r$ such that $V$ is the image of $H^m(X)_{\ell}(r) \to H^n(Y)_{\ell}(s)$.

(ii) The action of $\pi_1^{\log}(S)$ on $V$ factors through $\Gal(\overline k/k)$. 
\medskip

(1)  $H^m(X)(r)^{\flat}\mapsto $  the set of all morphisms $H^m(X)(r)^{\flat} \to H^n(Y)(s)^{\flat}$ such that the image of $H^m(X)_{\ell}(r)\to H^n(Y)_{\ell}(s)$ is contained in $V$. 

\end{para}

\begin{para} We expect that 
 this category (\MM)  is the category of mixed motives over $k$.

\end{para}
 
 \begin{para}
The above may be one of the simplest constructions of the category of mixed motives, and by the comparison with the Hodge version in \ref{Hodge}, we expect that the obtained category is the right one. 

However, it is not clear whether the above (\MM) contains the ``usual''  mixed motives $H^m(T)(r)$ associated with schemes $T$ of finite type over $k$. 
We give below another construction of the category of mixed motives over $k$ containing these ``usual'' objects, again by using log pure motives, and will conjecture that these two constructions give the same category. 

\end{para}

\begin{para}

For this, we use the homotopy $K$-theory $KH_n$ ($n\in \Z$, it is important for us that $n$ can be negative here) defined by Weibel \cite{W1}. There is a canonical homomorphism $K_n \to KH_n$ from Quillen's $K$-theory $K_n$ which is an isomorphism for regular Noetherian schemes. The reason why we use $KH_n$, not   Quillen's $K$-theory, is that we use the Riemann--Roch theorem for $KH_n$ proved in \cite{Na}. 

\end{para}

\begin{para}
For a scheme $T$ of finite type over $k$, let $H^m(T)_{\ell}=H^m_{\et}(T\otimes_k \overline k, \Q_{\ell})$. 

Let $X\in \cP$. Let $Y$ be an object of $\cP$  (resp.\  a scheme of finite type over $k$).  By a morphism $h: H^m(X)(r)^{\flat*}\to H^n(Y)(s)^{\flat*}$ (resp.\ $H^m(X)(r)^{\flat*}\to H^n(Y)(s)$)  of symbols, we mean a $\Q_{\ell}$-homomorphism $H^m(X)_{\ell}(r)\to H^n(Y)_{\ell}(s)$ obtained from some element of 
$$\gr^{d-r+s} KH_{(m-2r)-(n-2s),\lim}(Z)\otimes \Q,\quad\text{where}\;\;Z=X\times_S Y \;\;(\text{resp.}\; Z=X\times Y).$$
 Here $d$ is the dimension of $X$. Note that an element of this $K$-group goes by the Chern class map to $H^{2d-m+n}(Z)_{\ell}(d-r+s)$, and by K\"unneth formula and by Poincar\'e duality of $X$, to $\Hom_{\Q_{\ell}}(H^m(X)_{\ell}(r), H^n(Y)_{\ell}(s))$.

For such a morphism $h$ and for a morphism $g: H^{m(1)}(X_1)(r(1))^{\flat*} \to H^{m}(X)(r)^{\flat*}$ with $X_1\in \cP$, the composition $h\circ g: H^{m(1)}(X_1)_{\ell}(r(1)) \to H^n(Y)_{\ell}(s)$ is a morphism
$H^{m(1)}(X_1)(r(1))^{\flat*} \to H^n(Y)(s)^{\flat*}$ (resp.\ $H^{m(1)}(X_1)(r(1))^{\flat*} \to H^n(Y)(s)$). The identity map $ H^{m}(X)(r)^{\flat*}\to  H^{m}(X)(r)^{\flat*}$ is a morphism. These are proved in the same way as the non-$\flat$ case in Proposition \ref{comp}, by replacing the Riemann--Roch theorem for $K_0$ by the Riemann--Roch theorem for $KH_n$ in \cite{Na} which works for projective morphisms locally of complete intersection. 

Thus we have a category (\LM$\flat*$), and for a scheme $Y$ of finite type over $k$, we have a contra-variant functor $$H^n(Y)(s): H^m(X)(r)^{\flat*}\mapsto \{\text{morphisms}\; H^m(X)(r)^{\flat*}\to H^n(Y)(s)\}$$ from (\LM$\flat*$) to the category of $\Q$-vector spaces.

Let (\MM$*$) be the smallest full subcategory $\cC$ of the category of contra-variant functors  from (\LM$\flat*$) to the category of $\Q$-vector spaces satisfying the following conditions (i) and (ii).

(i) $\cC$ contains the functors
 $H^n(T)(s)$ for  schemes $T$ of finite type over $k$ and for $n,s\in \Z$.

(ii)  The kernel of every morphism of $\cC$ belongs to $\cC$.

\noindent 
That is, if $\cC_0$ denotes the category of  the functors $H^n(T)(s)$ for  schemes $T$ of finite type over $k$ and for $n,s\in \Z$ and if $\cC_{i+1}$ is the category of functors which are kernels of some morphisms of $\cC_i$, then  (\MM$*$)$=\bigcup_{i \geq 0}\; \cC_i$. 

Thus  (\MM$*$) is an additive category with  kernels of morphisms. The authors expect that it is an abelian category, but have not yet proved  it. The authors have not yet proved that the category $(\MM)$ is stable under taking kernels. 
\end{para}

\begin{para} For any scheme $T$ and for an integer $t\geq 0$, we have a canonical homomorphism $K_0(T\times {\Bbb G}_m^t)\to KH_{-t}(T)$, and the Chern class map  on the former $K$-group factors through the Chern class map on the latter $K$-group. Hence we have a functor
$$(\LM\flat)\to (\LM\flat*)$$
(the objects are the same but the set of morphisms might be enlarged  in the latter category).

\end{para} 

\begin{conj} 
$(\LM\flat)= (\LM\flat*)$ and $(\MM)= (\MM*)$.

\end{conj}

To check that our definitions of the category of mixed motives are reasonable, we show  an example \ref{sst} with our definitions for which the problems on Tate conjecture and Hodge conjecture  (\ref{Tate})  and the monodromy conjecture \ref{Nconj} on mixed motives have affirmative answers (Proposition \ref{p:Tate}).

\begin{para}\label{Tate} Let $Y$ and $Z$ be schemes of finite type over $k$ (resp.\ objects of $\cP$) and let $m,n,r,s\in \Z$. We ask whether the following (1) and (2) are true. 

$(1)$ (Tate conjecture.) Assume that $k$ is finitely generated over the prime field. 
Then
$$\Q_{\ell}\otimes_{\Q} \mathrm{Mor}_{(\MM*)}(H^m(Y)(r), H^n(Z)(s))\overset{\cong}\to \Hom_{\Gal(\overline k/k)}(H^m(Y)_{\ell}(r), H^n(Z)_{\ell}(s))$$ 
$$\text{(resp.}\ \Q_{\ell}\otimes_{\Q} \mathrm{Mor}_{(\LM\flat)}(H^m(Y)(r)^{\flat}, H^n(Z)(s)^{\flat})\overset{\cong}\to \Hom_{\pi_1^{\log}(S)}(H^m(Y)_{\ell}(r), H^n(Z)_{\ell}(s)){\text ).}$$ 

$(2)$ (Hodge conjecture.) Assume that $k=\C$. Then
$$\mathrm{Mor}_{(\MM*)}(H^m(Y)(r), H^n(Z)(s))\overset{\cong}\to \Hom_{(\MH)}(H^m(Y)(r)_H, H^n(Z)(s)_H)$$ 
$$\text{(resp.}\ \mathrm{Mor}_{(\LM\flat)}(H^m(Y)(r)^{\flat}, H^n(Z)(s)^{\flat})\overset{\cong}\to \Hom_{(\LH\flat)}(H^m(Y)(r)^{\flat}_H, H^n(Z)(s)^{\flat}_H))\text{).}$$ 
 Here $(\cdot)_H$ is the associated mixed Hodge structure (resp.\ mixed Hodge structure with $N$). 

\end{para}

\begin{para}\label{Bloch}

The conjectures in \ref{Tate} for the first isomorphisms in (1), (2) are in general false. The example  in Appendix of \cite{J} written by S.\ Bloch is a counter-example for the first isomorphism in (2) in which $Y=\Spec(\C)$, $m=0$, $r=0$, $Z$ is the  $W$ there which is three dimensional and singular, $n=4$, $s=2$. A counter-example for the first isomorphism in (1) is obtained from by defining this $W$ over a number field.

We expect that the conjectures for the second isomorphisms in (1), (2) are true in general. 
We expect that the above conjectures for the first isomorphisms in (1), (2) are true for smooth $Y$, $Z$, and more generally, for the underlying schemes of log smooth saturated schemes over the standard log point. 
\end{para}

\begin{rem}\label{Ja} (1) For  singular varieties, 
the Hodge conjecture [a Hodge class in homology $=$ an algebraic cycle class] and the Tate conjecture  [a Tate class (a Galois invariant element) in  homology $=$ an algebraic cycle class] are formulated in  Part II of Jannsen \cite{J} and are shown to be equivalent to the classical Hodge conjecture and Tate conjecture for projective smooth varieties (and hence are believed to be true), but the Hodge conjecture [a Hodge class in cohomology $=$ an algebraic cycle class] and the Tate conjecture [a Tate class in cohomology $=$ an algebraic cycle class] are false by Appendix of \cite{J}  written by Bloch. 
The  counter-examples in \ref{Bloch} appear because our theory considers cohomology $H^m(X)(r)$, not homology $H_m(X)(r)$. 

(2) In Part II of \cite{J}, for smooth varieties, conjectures [Hodge classes in cohomology come from Quillen's $K$-theory] and 
[Tate classes in cohomology come from   Quillen's $K$-theory] (for various Tate twists of the cohomology) are formulated. These are essentially 
the first isomorphisms in (1), (2) of \ref{Tate} for $Y=\Spec(k)$ and $m=0$, $r=0$, and $Z$ smooth, though we use the  homotopy $K$-theory $KH$, not Quillen's $K$-theory.

\end{rem}

\begin{conj}\label{Nconj} (Monodromy conjecture which tells that the monodromy operator comes from geometry, not only from Galois theory.)

   For $X\in \cP$, the monodromy operator $N: H^m(X)_{\ell}\to H^m(X)_{\ell}(-1)$ is a morphism $H^m(X)^{\flat}\to H^m(X)(-1)^{\flat}$, and hence is a morphism  $H^m(X)^{\flat*}\to H^m(X)(-1)^{\flat*}$.
 
\end{conj}

\begin{para}\label{sst} Example.  Let $\La$ be a discrete valuation ring with residue field $k$, and let $\frak X$ be a projective regular flat scheme over $\La$  of relative dimension one  with smooth generic fiber and with  semistable reduction. We assume that  the special fiber of $\frak X$ is a simple normal crossing divisor.
Endow $\Spec(\La)$ and $\frak X$ with the canonical log structures. We regard  $S$ as the closed point of $\Spec(\La)$  with the induced log structure. Let $X$ be the fs log scheme $\frak X \times_{\Spec(\La)} S$ over $S$. Then $X\in\cP$. Let $T$ be the underlying scheme $\frak X\otimes_{\La} k$ over $k$ of $X$.  
  We have a canonical injective homomorphism $H^1(T)_{\ell}\to H^1(X)_{\ell}$. 

\smallskip

\noindent Remark. $H^1(X)^{\flat}$ (or $H^1(X)^{\flat*}$) is regarded as the limit mixed motive, an analogue of limit mixed Hodge structure.

\end{para}

\begin{prop}
\label{p:Tate}
 Let the notation be as in {\rm \ref{sst}}. 

$(1)$ The Tate conjecture and the Hodge conjecture {\rm \ref{Tate}} for $(\MM*)$ are true in the case $Y=\Spec(k)$, $Z=T$, $m=0$, $n=1$, $r=s=0$.

$(2)$  The Tate conjecture and the Hodge conjecture {\rm \ref{Tate}} for $(\LM\flat)$ are true in the case $Y=S$, $Z=X$, $m=0$, $n=1$, $r=s=0$, and also in the case $Y=X$, $Z=S$, $m=1$, $n=0$, $r=0$, $s=-1$.

$(3)$ The monodromy operator $N: H^1(X)_{\ell}\to H^1(X)_{\ell}(-1)$ is a morphism of $(\LM\flat)$ (and hence a morphism of $(\LM\flat*)$). 

\end{prop}

\begin{pf} In (1) and (2), we only discuss the Tate conjecture. The proof for the Hodge conjecture is similar. 

In the discussion about $(\MM*)$ (resp.\ $(\LM\flat)$), we denote $H^0(\Spec(k))$ (resp.\ $H^0(S)^{\flat}$) by $\Q$. With this notation, the Tate conjecture in (1) is written as $$\Q_{\ell}\otimes_{\Q} \text{Mor}_{(\MM*)}(\Q, H^1(T))\overset{\cong}\to \Hom_{\Gal(\overline k/k)}(\Q_{\ell}, H^1(T)_{\ell}),$$ and the statements on the Tate conjecture in (2) are written as $$\Q_{\ell}\otimes_{\Q} \text{Mor}_{(\LM\flat)}(\Q, H^1(X)^{\flat})\overset{\cong}\to \Hom_{\pi_1^{\log}(S)}(\Q_{\ell}, H^1(X)_{\ell}),$$
$$\Q_{\ell}\otimes_{\Q} \text{Mor}_{(\LM\flat)}(H^1(X)^{\flat},\Q(-1))\overset{\cong}\to \Hom_{\pi_1^{\log}(S)}(H^1(X)_{\ell},\Q_{\ell}(-1)).$$

By Galois descent, we may and do assume that all singular points of $T$ are $k$-rational. Let $A$ be the set of all singular points of $T$ and let $B$ be the set of all generic points of $T$.  
 The following (i) and (ii) are well-known.

(i) We have  a canonical isomorphism $Q\overset{\cong}\to H^1_{\et}(T, \Z)$, where $Q$ is the cokernel of a natural homomorphism $\Z^B\to \Z^A$,  and it induces
an isomorphism from $\Q_{\ell} \otimes_{\Z} Q $ to the $G$-invariant part of $H^1(X)_{\ell}$, where $G=\pi_1^{\log}(S)$. Hence by  the Poincar\'e duality, we  have  an isomorphism from the $G$-coinvariant of $H^1(X)_{\ell}(1)$ to $\Q_{\ell}\otimes_{\Z} P$, where $P=\Hom(Q, \Z)$.

(ii) The monodromy logarithm $N: H^1(X)_{\ell}\to H^1(X)_{\ell}(-1)$ is the composition $${\rm (*)} \quad H^1(X)_{\ell}\to \Q_{\ell}(-1)\otimes_{\Z} P \to \Q_{\ell}(-1)^A \to \Q_{\ell}(-1)\otimes_{\Z} Q \to H^1(X)_{\ell}(-1).$$

By Theorem 3.3 and Theorem 5.1 of  \cite{W1} and by Lemma 2.3 of \cite{W2}, we have an isomorphism $KH_{-1}(T) \cong H^1_{\et}(T, \Z)$ and the Chern class map $KH_{-1}(T)\otimes \Q \to \gr^0KH_{-1}(T) \otimes \Q\to H^1(T)_{\ell}$ corresponds to the canonical map $H^1_{\et}(T, \Z)\to H^1(T)_{\ell}$. 
By definition, $\text{Mor}_{(\MM*)}(\Q, H^1(T))$   is the image of $KH_{-1}(T) \otimes \Q\to H^1(T)_{\ell}$ 
and hence we have (1). 

Next we consider (2). We prove first the version of the Tate conjecture in (2) in which we replace $(\LM\flat)$ by $(\LM\flat*)$. 
By definition, $\text{Mor}_{(\LM\flat*)}(\Q, H^1(X)^{\flat*})$ is the image of $KH_{-1}(T) \otimes \Q\to H^1(X)_{\ell}$ 
 and $\text{Mor}_{(\LM\flat*)}(H^1(X)^{\flat*}, \Q(-1))$  is the image of $KH_{-1}(T)\otimes\Q\to H^1(X)_{\ell}\cong \Hom(H^1(X)_{\ell}, \Q_{\ell}(-1))$ and hence we have the 
$(\LM\flat*)$ version of (2).

Now we consider $(\LM\flat)$. 
By definition, 
$\text{Mor}_{(\LM\flat)}(\Q, H^1(X)^{\flat})$ is the image of $\gr^1K_0(T\times {\Bbb G}_m) \otimes \Q\to H^1(X)_{\ell}$ and $\text{Mor}_{(\LM\flat)}(H^1(X), \Q(-1))$  is the image of $\gr^1K_0(T\times {\Bbb G}_m)\otimes\Q\to H^1(X)_{\ell}\cong \Hom(H^1(X)_{\ell}, \Q_{\ell}(-1))$. These Chern class maps factor through the above Chern class maps  on $\gr^1 KH_{-1}(T)\otimes \Q$. Since the map $\gr^1 K_0(T\times {\Bbb G}_m)\otimes \Q \to \gr^1 KH_{-1}(T)\otimes \Q \cong H^1_{\et}(T,\Z)\otimes \Q$  has a right inverse defined by 
$H^1_{\et}(T, \Z) \to H^1_{\et}(T\times {\Bbb  G}_m, {\Bbb G}_m) = \mathrm{Pic}\,(T\times {\Bbb G}_m)\to \gr^1K_0(T\times {\Bbb G}_m)$ in which the first arrow is the product with the coordinate function of ${\Bbb G}_m$,   we have  (2).

(3) follows from the above (ii) because every arrow in $(*)$ in (ii) is a morphism in (\MM). 
\end{pf}

\begin{para}
\label{**}
We have considered  the category of mixed motives modulo homological equivalence. The above method of the  construction of $(\MM*)$ works without homological equivalence as follows, by using the $K$-groups as the sets of morphisms. We define the modified version  $(\LM\flat\text{$**$})$ of $(\LM\flat*)$  as the category of symbols $h(X)(r)^{\flat\text{$**$}}$, where $X\in \cP$ and $r\in \Z$. We define the set of morphisms from $h(X)(r)^{\flat\text{$**$}}$ to $h(Y)(s)^{\flat\text{$**$}}$  to be $\bigoplus_{n\in \Z}\; \gr^{d-r+s}KH_{n, \lim}(X\times_S Y)\otimes \Q$,  where $d$ is the dimension of $X$. We define the modified version $(\MM\text{$**$})$ of $(\MM*)$ as 
 the category $\bigcup_{i\geq 0}\; \cC_i$ of contra-variant functors from $(\LM\flat\text{$**$})$ to the category of $\Q$-vector spaces, where $\cC_0$ is the category of the  functors 
$$h(T)(s):  
h(X)(r)^{\flat\text{$**$}}\mapsto \bigoplus_{n\in \Z}\;  \gr^{d-r+s}KH_{n, \lim}(X\times T) \otimes \Q$$  
for schemes $T$ of finite type over $k$ and for $s\in \Z$, and $\cC_{i+1}$ is the category of functors which are kernels of some  morphisms of $\cC_{i}$.

Thus (\MM$**$) is an additive category with kernels of morphisms. We expect that it is an abelian category.

\end{para}
\begin{para} In the case where the characteristic of $k$ is $0$, a definition of the category of mixed motives is given in Part I  of Jannsen \cite{J} by considering smooth (not necessarily proper) schemes. In his definition, a morphism of mixed motives is a compatible family of homomorphisms of various realizations (including the $\Q$-Betti realization; $K$-theory is not used in this definition). His definition and our definition are connected by the Tate conjecture for the first isomorphism in (1) of \ref{Tate} for smooth schemes. 

Our definition works also in positive characteristic in which we do not have the Betti realization.

The authors do not see how our definition is related to the work 
\cite{V} on mixed motives and 
\cite{BPO} on log motives. 

\end{para}


\begin{thebibliography}{99}

\bibitem{BPO} 
{\sc Binda, F., Park, D., \O stv{\ae}r, P.\ A.,}
{\em Triangulated categories of logarithmic motives over a field},
preprint, arXiv:2004.12298.   

\bibitem{B}
{\sc B\"uhler, T.}, {\em Exact categories}, 
Expositiones Math.\ {\bf 28} (2010), 1--69.

\bibitem{CKS1}
{\sc Cattani, E., Kaplan, A., Schmid, W.},
{\em Degeneration of Hodge structures},
Ann.\ of Math.\ {\bf 123} (1986), 457--535.

\bibitem{IKNU}{\sc Ito, T., Kato, K., Nakayama, C., Usui, S.},
 {\em On log motives},
Tunisian J.\ Math.\ {\bf 2} (2020), 733--789.  

\bibitem{J}
{\sc Jannsen, U.},
{\rm Mixed motives and algebraic $K$-theory},
Lecture Notes in Math. {\bf 1400}, Springer-Verlag, Berlin, New York, 1990.

\bibitem{Kas86}
{\sc Kashiwara, M.}, 
{\em A study of variation of mixed Hodge structure}, 
Publ.\ R.I.M.S., Kyoto Univ.\ {\bf 22} (1986), 991--1024.

\bibitem{KNU08}{\sc Kato, K., Nakayama, C., Usui, S.},
{\em $\SL(2)$-orbit theorem for degeneration of mixed Hodge structure},
 J.\ Algebraic Geometry {\bf 17} (2008), 401--479.

\bibitem{KNU3}
{\sc Kato, K., Nakayama, C., Usui, S.},
{\em Classifying spaces of degenerating mixed Hodge structures,
III$:$ Spaces of nilpotent orbits}, 
J.\ Algebraic Geometry {\bf 22} (2013), 671--772.

\bibitem{KNU4}
{\sc Kato, K., Nakayama, C., Usui, S.},
{\em Classifying spaces of degenerating mixed Hodge structures,
IV$:$ The fundamental diagram}, 
Kyoto J.\ Math.\ {\bf 58} (2018), 289--426.

\bibitem{KNU5}
{\sc Kato, K., Nakayama, C., Usui, S.},
{\em Classifying spaces of degenerating mixed Hodge structures,
V$:$ Extended period domains and algebraic groups},
preprint.

\bibitem{KU}
{\sc Kato, K., Usui, S.}, 
{\rm Classifying spaces of degenerating polarized 
Hodge structures}, 
Ann.\ Math.\ Studies {\bf 169}, Princeton Univ.\ Press, Princeton, NJ, 2009.

\bibitem{Nc}
{\sc Nakayama, C.}, 
{\em Logarithmic \'etale cohomology}, 
Math.\ Ann.\ {\bf 308} (1997), 365--404.

\bibitem{Nc2}
{\sc Nakayama, C.}, 
{\em Logarithmic \'etale cohomology, II}, 
Advances in Math.\ {\bf 314} (2017), 663--725.

\bibitem{Na}
{\sc Navarro, A.}, 
{\em Riemann--Roch for homotopy invariant $K$-theory and Gysin morphisms},
Advances in Math.\ {\bf 328} (2018), 501--554. 

\bibitem{P}
{\sc Positselski, L.}, 
{\em Mixed Artin--Tate motives with finite coefficients}, Moscow Math.\ J.\ 
{\bf 11}  (2011).

\bibitem{Sc}
{\sc Schmid, W.},
{\em  Variation of Hodge structure{\rm:}
the singularities of the period mapping},
 Invent.\ Math.\ {\bf 22} (1973), 
211--319.

\bibitem{V}
{\sc Voevodsky, V.}, 
{\em Triangulated categories of motives over a field},  
In: Cycles, transfers, and motivic homology theories, 
Ann.\ 
of Math.\ Studies {\bf 143}, Princeton University Press, Princeton, NJ, 2000, 188--238.

\bibitem{W1}
{\sc Weibel, C. A.}, 
{\em Homotopy algebraic $K$-theory},
Contemporary Math.\ {\bf 83} (1989), 461--488. 

\bibitem{W2}
{\sc Weibel, C. A.}, 
{\em The negative K-theory of normal surfaces}, Duke Math.\ J.\ {\bf 108} (2001), 1--35.
\end{thebibliography}
\end{document}